\documentclass[a4paper]{amsart}
\usepackage{amsfonts,amssymb,mathrsfs,bbold,comment}
\usepackage{stmaryrd,soul,adjustbox,circuitikz}
\usepackage[a4paper,hmargin=26mm,vmargin=28mm]{geometry}
\usepackage{enumitem}
\setlist[enumerate]{label=\upshape(\alph*)., ref=\alph*}
\usepackage[svgnames]{xcolor}
\usepackage{booktabs,upgreek}
\usepackage{xparse,mathtools,etoolbox}
\usepackage{atableau}
\usepackage{tikz}
\usetikzlibrary{cd,arrows,decorations.markings}
\usepackage{todonotes}

\title{A Specht Filtration of Permutation Modules Over KLR Algebras}
\author{Tao Qin}
\address[T.Q.]{The University of Sydney, School of Mathematics and Statistics F07, NSW 2006, Australia}
\email{tao.qin97@gmail.com}


\usepackage{cite}
\usepackage{hyperref}
\hypersetup{
  pdftitle={Specht Filtration of Permutation Module}, 
  pdfauthor={Tao Qin},
  colorlinks=true,
  linkcolor=blue,
  anchorcolor=red,
  citecolor=blue,
  urlcolor=blue
}

\newcommand\enumref[2]{\hyperref[#2]{\autoref*{#1}(\ref*{#2})}}
\let\<=\langle
\let\>=\rangle
\def\({\big(}
\def\){\big)}

\newcommand{\Sym}{\mathfrak S}

\newcommand{\Z}{\mathbb{Z}}

\DeclarePairedDelimiterX\cart[2]{\langle}{\rangle}{#1,#2}

\newcommand{\Par}[1][\Lambda]{\mathscr{P}^{#1}}

\newcommand\blam{{\boldsymbol\lambda}}
\newcommand\bmu{{\boldsymbol\mu}}

\def\c{\mathsf{c}}

\newcommand\bs{\mathbf{s}}
\newcommand\bk{\mathbb{k}}

\newcommand\height{\text{ht}}

\DeclarePairedDelimiterX{\set}[1]{\{}{\}}{\setargs{#1}}
\NewDocumentCommand{\setargs}{>{\SplitArgument{1}{|}}m}{\setargsaux#1}
\NewDocumentCommand{\setargsaux}{mm}
{\IfNoValueTF{#2}{#1} {#1\,\delimsize|\,\mathopen{}#2}}%{#1\:;\:#2}

\def\pmod#1{\text{ }(\textrm{mod } #1)\,}

\usepackage{aliascnt}
\def\NewTheorem#1{%
  \newaliascnt{#1}{equation}%
  \newtheorem{#1}[#1]{#1}%
  \aliascntresetthe{#1}%
  \expandafter\def\csname #1autorefname\endcsname{#1}%
}
\def\equationautorefname~#1\null{(#1)\null}
\def\itemautorefname~#1\null{(#1)\null}

\numberwithin{equation}{section}
\NewTheorem{Proposition}
\NewTheorem{Theorem}
\NewTheorem{Conjecture}
\NewTheorem{Condition}
\NewTheorem{Corollary}
\NewTheorem{Lemma}
\NewTheorem{Definition}
\theoremstyle{definition}
\NewTheorem{Notation}
\NewTheorem{Example}
\AtEndEnvironment{Example}{\null\hfill$\Diamond$}%
\NewTheorem{Examples}
\AtEndEnvironment{Examples}{\null\hfill$\Diamond$}%
\theoremstyle{remark}
\NewTheorem{Remark}
\AtEndEnvironment{Remark}{\null\hfill$\Diamond$}%
\NewTheorem{Claim}

\newcommand\bi{\mathbf{i}}
\newcommand\bj{\mathbf{j}}

\newcommand{\DeclareMyOperator}[1]{%
  \expandafter\DeclareMathOperator\csname #1\endcsname{#1}
}
\forcsvlist{\DeclareMyOperator}{%
   Mat,Sh,Shaded,cont,Add,
  ch,
  diag,
  End,
  END,
  head,
  Hom,
  HOM,
  im,
  inv,
  Irr,
  Deg,
  Ind,
  Rad,
  Rem,
  rad,
  Res,
  res,
  soc,
  supp,
  Shape,
  Std,
  ST,
  RST,
  CST,
  Top,
}

\newcommand\ShT[1][l]{\mathop{ShT}}

\DeclareMathOperator\arrow{\text{---}}
\let\kill\relax
\DeclarePairedDelimiterX\kill[2]{(}{)}{#1|#2}

\usepackage{mathrsfs}

\newcommand\Rn[1][n]{\mathscr{R}^\Lambda_{#1}}

\NewDocumentCommand\Rnaf{ D<>{n}}{\Rn[]<#1>}
\DeclareRobustCommand{\Aone}[1][e+1]{\ensuremath{A^{(1)}_{#1}}}
\newcommand\Cone[1][e]{C^{(1)}_{#1}}
\newcommand\Atwo[1][2e]{A^{(2)}_{#1}}
\newcommand\Dtwo[1][e+1]{D^{(2)}_{#1}}

\newcommand{\al}{\alpha}

\renewcommand{\phi}{\varphi}

\tikzset{
  centered/.style = {
     baseline = {([yshift=#1]current bounding box.center)}
  },
  centered/.default={-0.5ex},
  ->-/.style={
    decoration={
      markings,
      mark=at position 0.6 with {\arrow[thin,scale=2]{>}}
    },
    postaction={decorate}
  },
  -<-/.style={
    decoration={
      markings,
      mark=at position 0.6 with {\arrow[thin,scale=2]{<}}
    },
    postaction={decorate}
  },
  circled/.style = {fill=Tan},
  domstyle/.style = {
    thick,
    scale=#1,
  },
  vertex/.style = {
    circle,
    ball color=MidnightBlue,
    font=\small,
    inner sep=0pt,
    minimum size=2mm
  },
  pics/domm/.style = {
    code = {
      \draw[domstyle=#1](0.8,0)--++(0.6,0.4)--++(0.6,-0.4)--++(-0.6,-0.4)--++(-0.6,0.4);
    }
  },
  pics/dom/.style = {
    code = {
      \draw[domstyle=#1](0.2,0)--(0.8,0)
        --++(0.6,0.4)--++(0.6,-0.4)--++(-0.6,-0.4)--++(-0.6,0.4);
    }
  },
  pics/domeq/.style = {
    code = {
      \pic at (0,0){dom=#1};
      \draw[scale=#1](0.8,-0.5)--++(1.2,0);
    }
  },
  pics/domneq/.style = {
    code = {
      \pic at (0,0){domeq=#1};
      \draw[scale=#1](0.9,-0.75)--++(1.0,0.5);
    }
  },
  pics/mod/.style = {
    code = {
      \draw[domstyle=#1](0,0)--++(-0.6,0)--++(-0.6,-0.4)--++(-0.6,0.4)--++(0.6,0.4)--++(0.6,-0.4);
    }
  },
  pics/Sdom/.style = {
    code = {
      \draw[domstyle=#1](0.2,0)--(0.8,0);
      \draw[scale=#1,fill=black](0.8,0)--++(0.6,0.4)--++(0.6,-0.4)--++(-0.6,-0.4)--++(-0.6,0.4);
    }
  },
  pics/Sdomeq/.style = {
    code = {
      \pic at (0,0){Sdom=#1};
      \draw[scale=#1](0.8,-0.5)--++(1.2,0);
    }
  },
  pics/Sdomneq/.style = {
    code = {
      \pic at (0,0){Sdomeq=#1};
      \draw[scale=#1](0.9,-0.75)--++(1.0,0.5);
    }
  },
}

\usepackage{enumitem}
\setlist[enumerate]{label=\upshape\alph*), ref=\alph*}
\let\realItem\item
\NewDocumentCommand\centeredItem{so}{% * disables hfil
   \IfNoValueTF{#2}{\realItem}{\realItem[#2]}%
   \IfBooleanF{#1}{\hfil}%
}
\newlist{Enumerate}{enumerate}{1}
\setlist[Enumerate]{% enumerate with references tied to equation number
  label=\textup{\alph*)},
  ref=\theequation\alph*,
}
\newlist{relations}{enumerate}{1}
\setlist[relations]{% only suited to "one line" items because of centering
  label=$(\text{KLR}_{\arabic*})$,
  before=\let\item\centeredItem
}

\begin{document}

\begin{abstract}
    Kleshchev–Mathas–Ram give a presentation of the Specht module $S^\lambda$ as a quotient of the permutation module $M^\lambda$. In this paper, we construct a (graded) Specht filtration of the permutation module $M^\lambda$ in the following cases: when $\lambda = (k,1^r)$ is a hook partition, over the KLR algebra of type $\Aone[e-1]$ for $e > 2$; and when $\lambda = (k,r)$ is a two-row partition with $k \geq r$, over the KLR algebra of type $A_\infty$. Furthermore, when $\lambda$ is an arbitrary partition in type $A_\infty$, we construct a filtration of $M^\lambda$ such that each subquotient $M_i / M_{i+1}$ admits a Specht resolution.
\end{abstract}
\maketitle

\setcounter{tocdepth}{1}
\tableofcontents
%%%%%%%%%%%%%%%%%%%%%%%%%%%%%%%%%%%%%%%%%%%%%%%%%%%%%%%%%%%%%%%%%%%%%%%%%%%%%%%%%%%%%%%%%%%%%%%%%%%%%%%%%%%%%%%%%%%%%%%%%%%%%%%%%%%%%%%%%%%%%%%%%%%%%%%%%%%%%%%%%%%%%%%%%%%%%%%%

\section{Introduction}
Khovanov, Lauda, and Rouquier independently introduced KLR algebras to categorify quantum groups; see \cite{khovanovlauda-klr-1, khovanovlauda-KLR-3} and \cite{rouquier-2kacmoody}. Several years later, Kang and Kashiwara \cite{kangkashiwara-klr-categorification} proved that the cyclotomic KLR algebras of any symmetrizable type categorify the highest weight modules of the corresponding quantum groups.

There are several well-known results concerning KLR algebras. For instance, the Brundan–Kleshchev isomorphism \cite{bk-blocks-iso} establishes an equivalence between Ariki–Koike algebras and cyclotomic KLR algebras of type $\Aone[e-1]$. Additionally, Hu and Mathas \cite{humathas-graded-cellular} constructed a graded cellular basis for cyclotomic KLR algebras of type $\Aone[e-1]$, and showed that the corresponding graded cell modules are isomorphic to the graded Specht modules constructed by Brundan, Kleshchev, and Wang \cite{bkw-graded-specht}. Furthermore, in \cite{hushi-dimension-klr}, Hu and Shi derived a dimension formula for cyclotomic KLR algebras of arbitrary symmetrizable type.

More recently, Evseev and Mathas \cite{evseevmathas-klr-deformation} introduced a deformation method to construct cellular bases for KLR algebras of both type $\Aone[e-1]$ and type $\Cone[e]$. It is also worth noting that Bowman \cite{bowman-many-cellular} used idempotent truncations of KLRW algebras (introduced by Webster \cite{webster-weighted-klr} as generalizations of KLR algebras) to construct a family of diagrammatic cellular bases for KLR algebras of type $\Aone[e-1]$. Mathas and Tubbenhauer have further explored the (sandwiched) cellularity of KLRW algebras in types $\Cone[e]$, $\Atwo[]$, and $\Dtwo[]$ \cite{mathastubbenhauer-klrw-ac, mathastubbenhauer-klrw-bad}, as well as in finite types \cite{mathastubbenhauer-klrw-finite}.

In \cite{kmr-universal-specht-type-A}, Kleshchev, Mathas, and Ram gave a highest-weight presentation of the graded Specht module, referred to as the \emph{universal Specht module}. The Specht module $S^\lambda$ is constructed as the quotient of the permutation module $M^\lambda$ by the ideal generated by the Garnir relations. Mathas conjectured that the permutation module $M^\lambda$ admits a filtration by Specht modules; equivalently, that the ideal generated by the Garnir relations admits such a Specht filtration. In this paper, we prove this conjecture for hook partitions $\lambda = (k, 1^r)$ in type $\Aone[e-1]$ with $e > 2$ (see \autoref{thm-specht-filtration-hook}). We also prove the conjecture in type $A_\infty$ for two-row partitions (see \autoref{thm:specht-filtration-2row}) and for arbitrary partitions $\lambda$, with appropriate modifications as stated in \autoref{thm:specht-filtration-general-partition}.

It is worth noting that the permutation module $M^{\lambda}$ can be identified with a higher-level Specht module by separating the rows of the partition $\lambda$ and regarding them as components of the multi-partition; see the discussion at the beginning of \autoref{sec:filtration}. Therefore, the (generalized) Specht filtration constructed in this paper can also be viewed as a (generalized) Specht filtration of a higher-level Specht module.

Several works have followed \cite{kmr-universal-specht-type-A}. Muth \cite{muth-graded-skew-specht} extended their construction of Specht modules to skew partitions and demonstrated that, in affine type A, cuspidal modules associated with a balanced convex preorder are skew Specht modules for certain hook shapes, in the context of the theory of cuspidal systems; see \cite{kleshchev-cuspidal-affine-klr,km-imaginary-schur-weyl,mcnamara-cuspidal-symmetric-affine} for classical results and \cite{admpss-cuspidal-ribbon-tableaux-type-a,mnss-skew-rock-cuspidal} for recent developments. Loubert \cite{loubert-homo-specht-hook} investigated homomorphisms between Specht modules $S^\lambda$ and $S^\mu$, where $\mu$ is a hook partition, for KLR algebras of type $\Aone[e-1]$ with $e > 2$. Subsequently, Hudak \cite{hudak-homo-specht-hook-e2} addressed the case $e = 2$. In \cite{aps-universal-specht-type-c}, Ariki, Park, and Speyer introduced the universal Specht module for type $C_\infty$ and proposed a conjecture in type $\Cone[e]$.  Twisted types such as $\Atwo[2e]$ and $\Dtwo[e+1]$ are also under active investigation; indeed, the motivation for this paper partially arose from efforts to understand the structure of universal Specht modules in type $\Atwo[2e]$.

The paper is organized as follows. In \autoref{sec:preliminary}, we introduce the necessary combinatorics and Lie theory. In \autoref{sec:klr-specht}, we describe the KLR algebras and the universal Specht modules. Experienced readers may wish to proceed directly to \autoref{sec:filtration}, where we present the construction of the Specht filtration (\autoref{thm-specht-filtration-hook}) for hook partitions in type~$\Aone[e-1]$ with $e > 1$. In \autoref{sec:specht-filtration-2row}, we construct a Specht filtration (\autoref{thm:specht-filtration-2row}) for two-row partitions in type~$A_\infty$. Furthermore, in \autoref{sec:specht-filtration-general-partition}, we develop a generalized Specht filtration (\autoref{thm:specht-filtration-general-partition}) of $M^\lambda$ for arbitrary partitions in type~$A_\infty$, extending the work from \autoref{sec:specht-filtration-2row}. Finally, in \autoref{sec:higher-level-skew}, we briefly discuss the higher-level case and introduce a skew Specht filtration. 
%%%%%%%%%%%%%%%%%%%%%%%%%%%%%%%%%%%%%%%%%%%%%%%%%%%%%%%%%%%%%%%%%%%%%%%%%%%%%%%%%%%%%%%%%%%%%%%%%%%%%%%%%%%%%%%%%%%%%%%%%%%%%%%%%%%%%%%%%%%%%%%%%%%%%%%%%%%%%%%%%%%%%%%%%%%%%%%%%%%%%%%%%%%%%%%%%%%%%%%%%%%%%%%%%%%%%%%%%%%%%%%%%%%%%%%%%%%%%%%%%%%%%%%%%%%%%%%%%%%%%%%%%%%%%%%%%%%%%%%%%%%%%%%%%%%%%%%%%%%%%%%%%%%%%%%%%%%%%%%%%%%%%%%%%%%%%%%%%%%%%%%%%%%%%%%%%%%%%%%%%%%%%%%%%%%%%%%%%%%%%%%%%%%%%%%%%%%%%%%%%%%%%%%%%%%%%%%%%%%%%%%%%%%%%%%%%%%%%%%%%%%%%%%%%%%%%%%%%%%%%%%%%%%%%%%%%%%%%%%%%%%%%%

\noindent\textbf{Acknowledgments.}
The author thanks the referees for their careful reading of the manuscript and for many helpful comments and suggestions.

We are grateful to Andrew Mathas for proposing the problem and for numerous insightful discussions; the Young-diagram package used in this work is also due to him. We also thank Huang Lin for valuable conversations on universal Specht modules in other types, and Nick Bridger, Shixuan Wang, Finn Klein, Tom Goertzen, Tasman Fell, and Joe Newton for helpful discussions. We also thank Liron Speyer and Chris Bowman for their careful reading and constructive feedback on this manuscript, which forms part of a PhD thesis. This work was partially supported by the Australian Research Council (ARC) Discovery Grant DP240101809.
%%%%%%%%%%%%%%%%%%%%%%%%%%%%%%%%%%%%%%%%%%%%%%%%%%%%%%%%%%%%%%%%%%%%%%%%%%%%%%%%%%%%%%%%%%%%%%%%%%%%%%%%%%%%%%%%%%%%%%%%%%%%%%%%%%%%%%%%%%%%%%%%%%%%%%%%%%%%%%%%%%%%%%%%%%%%%%%%%%%%%%%%%%%%%%%%%%%%%%%%%%%%%%%%%%%%%%%%%%%%%%%%%%%%%%%%%%%%%%%%%%%%%%%%%%%%%%%%%%%%%%%%%%%%%%%%%%%%%%%%%%%%%%%%%%%%%%%%%%%%%%%%%%%%%%%%%%%%%%%%%%%%%%%%%%%%%%%%%%%%%%%%%%%%%%%%%%%%%%%%%%%%%%%%%%%%%%%%%%%%%%%%%%%%%%%%%%%%%%%%%%%%%%%%%%%%%%%%%%%%%%%%%%%%%%%%%%%%%%%%%%%%%%%%%%%%%%%%%%%%%%%%%%%%%%%%%%%%%%%%%%%%%%%%%%%%%%%%%%%%%%%%%%%%%%%%%%%%%%%%%%%%%%%%%%%%%%%%%%%%%%%%%%%%%%%%%%%%%%%%%%%%%%%%%%%%%%%%%%%%%%%%%%%%%%%%%%%%%%%%%%%%%%%%%%%%%%%%%%%%%%%%%%%%%%%%%%%%%%%%%%%%%%%%%%%%%%%%%%%%%%
\section{Preliminaries}\label{sec:preliminary}
We collect some notation that will be used throughout this paper. For further details, the interested reader is referred to the papers cited in the introduction. In particular, all notation used here can be found in \cite{kmr-universal-specht-type-A}.
\subsection{Cartan data}\label{subsec:cartan-data}
In this paper, we consider the quivers $\Aone[e-1], (e>2)$:
\begin{center}
    \begin{tikzpicture}
        % Draw edges
        \draw(2.5,0)--(4,0)--(2,0.5)--(0,0)--(1.5,0);
        \draw[dashed](1.5,0)--(2.5,0);

         % Draw and label nodes
         \shade[ball color=black] (0,0) circle(3pt);   % label 1
         \node[below] at (0,0) {$1$};

        \shade[ball color=black] (1,0) circle(3pt);   % label 2
         \node[below] at (1,0) {$2$};

        \shade[ball color=black] (2,0.5) circle(3pt); % label 0 (top node)
        \node[above] at (2,0.5) {$0$};

         \shade[ball color=black] (3,0) circle(3pt);   % label e-1
        \node[below] at (3,0) {$e\!-\!2$};

         \shade[ball color=black] (4,0) circle(3pt);   % label e
        \node[below] at (4,0) {$e\!-\!1$};
    \end{tikzpicture}
\end{center}

 and $A_\infty$:
\begin{center}
    \begin{tikzpicture}
        % Draw edges
        \draw[dashed](0,0)--(-1,0);
         \draw(0,0)--(1,0);
        \draw(1,0)--(1.5,0);
        \draw[dashed](1.5,0)--(2.5,0);
        \draw(2.5,0)--(3,0);
        \draw(3,0)--(4,0);
        \draw[dashed](4,0)--(5,0);

        % Draw and label nodes: positions and corresponding labels
        \shade[ball color=black](0,0) circle(3pt);   % label 0
        \node[below] at (0,0) {$-1$};
    
        \shade[ball color=black](1,0) circle(3pt);   % label -1
        \node[below] at (1,0) {$0$};

        \shade[ball color=black](3,0) circle(3pt);   % label n
        \node[below] at (3,0) {$n$};

        \shade[ball color=black](4,0) circle(3pt);   % label n+1
        \node[below] at (4,0) {$n\!+\!1$};
    \end{tikzpicture}
\end{center}

Since the $A_\infty$ case in this paper is just a trivial subcase of the $\Aone[e-1]$ case, we may focus on the case $\Aone[e-1]$. 

Let $\Gamma$ be the quiver of type $\Aone[e-1]$. It has vertex set $I=\{0,1,\cdots,e-1\}$ which can be identified with $\mathbb{Z}/e\mathbb{Z}$. In particular, we identify $e$ with $0$.

Throughout this paper, we always fix an orientation by $i\rightarrow i+1$ for each $i$. 

The \emph{affine Cartan matrix} $(a_{ij})_{i,j\in I}$ of type $\Aone[e-1]$ is the following $e\times e$-matrix:
\[
\begin{pmatrix}\label{cartan-matrix-aone}
  2 & -1 & 0 & \cdots & 0 & 0 & -1 \\
  -1 & 2 & -1 & 0 & \cdots & 0 & 0 \\
  0 & -1 & 2 & -1 & \cdots & 0 & 0 \\
  \vdots & \vdots & \vdots & \ddots & \ddots & \vdots & \vdots \\
  0 & \cdots & 0 & -1 & 2 & -1 & 0 \\
  0 & \cdots & 0 & 0 & -1 & 2 & -1 \\
  -1 & 0 & \cdots & 0 & 0 & -1 & 2 \\
\end{pmatrix}
\]

The \emph{simple roots} are $\{\al_i\mid i\in I\}$ and $Q^+ := \bigoplus_{i \in I} \Z_{\geq 0} \alpha_i$ is the positive part of the root lattice. For $\alpha \in Q^+$ let $\height(\alpha)$ be the {\em height of~$\al$}. That is, $\height(\al)$ is the sum of the coefficients when $\al$ is expanded in terms of the $\alpha_i$'s.

Let~$\Sym_n$ be the symmetric group on~$n$ letters and let $\sigma_r = (r, r+1)$, for $1\leq r < n$, be the simple transpositions of~$\Sym_n$. Then~$\Sym_n$ acts on the left on the set~$I^n$ by place permutations.
If $\bi = (\bi_1, \dots , \bi_n) \in I^n$ then its {\em weight} is
$|\bi| := \alpha_{i_1} + \cdots + \alpha_{i_n} \in Q^+$.  The $\Sym_n$-orbits on $I^n$ are the sets 
\[
  I^\alpha := \{\bi \in I^n\mid \al=|\bi|\} 
\]
parametrized by all $\alpha \in Q^+$ of height $n$.

Let $\Lambda_0,\cdots,\Lambda_{e-1}$ be the \emph{fundamental weights}. The \emph{dominant weight lattice} is then \[P^+:=\bigoplus_{i \in I} \Z_{\geq 0} \Lambda_i\]
Any element $\Lambda=\sum\limits_{0\leq i\leq e-1}a_i\Lambda_i\in P^+$ is a \emph{dominant weight} and the sum of coefficients $\sum\limits_{0\leq i\leq e-1}a_i$ is the \emph{level} of $\Lambda$. 
%%%%%%%%%%%%%%%%%%%%%%%%%%%%%%%%%
\subsection{Partitions and tableaux}
In this section, we fix a quiver $\Aone[e-1]$ with vertex set $I$. Let $\Par[]_n$ be the set of all partitions of $n$ and put $\Par[]:=\bigsqcup_{n\geq 0}\Par[]_n$. 
The {\em Young diagram} of the partition $\lambda=(\lambda_1,\cdots,\lambda_r)\in \Par[]$ is 
\[
[\lambda]:=\{(a,b)\in\Z_{>0}\times\Z_{>0}\mid 1\leq b\leq \lambda_a,1\leq a\leq r\}.
\]
The elements of this set
are the {\em nodes of $[\lambda]$}. Let $|\lambda|=\sum\limits_{1\leq i\leq r}\lambda_r$ be the size of $\lambda$. 

Fix a fundamental weight $\Lambda:=\Lambda_j$ for some $j\in I$. To each node $A=(a,b)$ we associate its \emph{residue}, which is the following element of $I=\Z/e\Z$:
\[
    \res_{\Lambda} A=(b-a+j)\pmod{e}.
\]
If there is no ambiguity, we may omit the weight $\Lambda$ and write $\res A$ instead.

An {\em $i$-node} is a node of residue $i$.
Define the {\em residue content of $\lambda$} to be
\[
    \alpha_\lambda:=\sum_{A\in [\lambda]}\al_{\res_{\Lambda} A} \in Q^+.
\]

Let $\Par[\Lambda]_n$ be the set of partitions of $n$ with the residue function $\res_{\Lambda}$ defined for each node of its Young diagram. Let $\Par[\Lambda]:=\bigsqcup_{n\geq 0}\Par[\Lambda]_n$. For each $\alpha\in Q^{+}$, set
\[
    \Par[\Lambda]_\alpha:=\{\lambda\in\Par\mid \alpha_{\lambda}=\al\}.
\]

An $\ell$-partition $\blam=(\blam^{(1)}\mid\cdots\mid\blam^{(\ell)})$ of $n$ is an $\ell$-tuple of partitions satisfying $\sum_{i=1}^{\ell}|\blam^{(i)}|=n$. Here $\blam^{(i)}$ is the \emph{$i$th component}, and the Young diagram $[\blam]$ is the following set:
\[
    \{(m,r,c)\mid 1\leq m\leq \ell, (r,c)\in [\blam^{(i)}]\}
\]

Given a dominant weight $\Lambda\in P^+$ of level $\ell$, an \emph{(integral) charge of $\Lambda$} is an $\ell$-tuple $\kappa=(\kappa_1,\dots,\kappa_\ell)\in\mathbb{Z}^\ell$ such that 
\[
    \sum_{i=1}^{\ell}\Lambda_{\overline{\kappa_i}}=\Lambda,\qquad \kappa_i\equiv\overline{\kappa_i}\pmod{e}
\]
Define the associated residue function $\res_{\kappa}$ on nodes of $[\blam]$ as follows. For a node $A=(m,r,c)\in[\blam]$ (so that $(r,c)\in[\blam^{(m)}]$), set
\[
    \res_{\kappa} A \;:=\; \res_{\Lambda_{\overline{\kappa_m}}}(r,c),
\]
Let $\alpha_{\blam^{(m)}}$ be the residue content of each component $\blam^{(m)}$, and define the residue content of $\blam$ to be $\alpha_{\blam}=\sum_{m=1}^{\ell}\alpha_{\blam^{(m)}}$. The residue of the first node of $\blam^{(m)}$ is called the \emph{leading residue} of $\blam^{(m)}$.

Similarly to the partition case, let $\Par[\kappa]_n$ be the set of $\ell$-partitions with the residue function on nodes determined by $\kappa$ as above. Let $\Par[\kappa]:=\bigsqcup_{n\geq 0}\Par[\kappa]_n$. For $\alpha\in Q^+$, set $\Par[\kappa]_{\alpha}=\{\blam\in \Par[\kappa]\mid \alpha_{\blam}=\alpha\}$.

In this paper, when the order of the components of an $\ell$-partition is fixed and the leading residues are already determined as $(i_1, \cdots, i_\ell)\in I^\ell$, we always fix a charge $\kappa$ satisfying $\kappa_1 > \kappa_2 > \cdots > \kappa_\ell$ and $\overline{\kappa_j} \equiv i_j \pmod{e}$ for each $j$. By abuse of notation, we often write $\Par[\Lambda]_{\alpha}:=\Par[\kappa]_{\alpha}$ in this case.

\medskip
Fix $\Lambda\in P^+$ and a charge $\kappa$ of $\Lambda$, and take $\blam\in\Par[\Lambda]_n$.
A \emph{$\blam$-tableau} $T$ is a bijection from the nodes of $[\blam]$ to the set $\{1,2,\cdots,n\}$. Informally, it is obtained from $[\blam]$ by inserting the integers $1,\dots,n$ into the nodes, with no repeats.
If the node $A=(m,r,c)\in[\blam]$ is occupied by the integer $k$ in $T$, then we write $T(m,r,c)=k$ and set $\res_T(k)=\res_{\kappa} A$. The \emph{residue sequence} of~$T$ is
\[
    \res(T)=\bi(T)=(i_1,\dots,i_n)\in I^n,
\]
where $i_k=\res_T(k)$ is the residue of the node occupied by $k$ in $T$ ($1\leq k\leq n$).

A $\blam$-tableau $T$ is \emph{row-standard} (resp.\ \emph{column-standard}) if its entries increase from left to right (resp.\ from top to bottom) along the rows (resp.\ columns) of each component of $T$.
A $\blam$-tableau $T$ is \emph{standard} if it is both row- and column-standard.
Let $\Std(\blam)$ be the set of standard $\blam$-tableaux.
  
Given $\blam\in\Par[\kappa]_n$ and $T\in\Std(\blam)$, the \emph{degree} of $T$ is defined inductively in \cite[Section~3.5]{bkw-graded-specht}. Since gradings play a minimal role in this paper, we omit further discussion of them.

The group $\Sym_n$ acts on the set of $\blam$-tableaux from the left by acting on the entries of the tableaux. Here we follow the convention that elements of $\Sym_n$ compose from right to left.

Let $T^\blam$ be the $\blam$-tableau in which the numbers $1,2,\dots,n$ appear in order from left to right along the successive rows, working from the top row to the bottom row and from the first component $\blam^{(1)}$ to the last component $\blam^{(\ell)}$. Set 
\begin{equation}\label{eq:residue-sequence}
    \bi^{\blam}:=\bi(T^{\blam}).
\end{equation} 
In \autoref{sec:filtration}, we often set $\bi=\bi^{\blam}$ if the $\ell$-partition $\blam$ is fixed.

For each row-standard $\blam$-tableau $T$ define the permutation $w^T \in \Sym_n$ by the equation
\begin{equation}\label{EWT}
w^T  T^\lambda=T.
\end{equation}
\begin{Example}\label{eg:tableau}
    Let $\lambda=(4,2)$, then the Young diagram $[\lambda]$ is:
    \[
        \Tableau{~~~~,~~}
    \]
    Fix quiver $\Aone[2]$ and let $\Lambda=\Lambda_0$, we fill the Young diagram by its residue:
    \[
        \Tableau{0120,20}
    \]
    The $\lambda$-tableau $T^\lambda$ is the following:
    \[
        \Tableau{1234,56}
    \]
    and $\bi^\lambda=\res(T^\lambda)=(0,1,2,0,2,0)$. Let $\sigma_4\in\Sym_6$ be the transposition $(4,5)$, then $\sigma_4 T^\lambda$ is the following tableau:
    \[
        \Tableau{1235,46}
    \]
    and $\res(\sigma_4 T^\lambda)=(0,1,2,2,0,0)$.
    Let $T$ be the following $\lambda$-tableau:
    \[
        \Tableau{1345,26}
    \]
    Then $w^T=\sigma_2\sigma_3\sigma_4$.
\end{Example}
%%%%%%%%%%%%%%%%%%%%%%%%%%%%%%%%%%%%%%%%%%%%%%%%%%%%%%%%%%%%%%%%%%%%%%%%%%%%%%%%%%%%%%%%%%%%%%%%%%%%%%%%%%%%%%%%%%%%%%%%%%%%%%%%%%%%%%%%%%%%%%%%%%%%%%%%%%%%%%%%%%%%%%%%%%%%%%%%%%%%%%%%%%%%%%%%%%%%%%%%%%%%%%%%%%%%%%%%%%%%%%%%%%%%%%%%%%%%%%%%%%%%%%%%%%%%%%%%%%%%%%%%%%%%%%%%%%%%%%%%%%%%%%%%%%%%%%%%%%%%%%%%%%%%%%%%%%%%%%%%%%%%%%%%%%%%%%%%%%%%%%%%%%%%%%%%%%%%%%%%%%%%%%%%%%%%%%%%%%%%%%%%%%%%%%%%%%%%%%%%%%%%%%%%%%%%%%%%%%%%%%%%%%%%%%%%%%%%%%%%%%%%%%%%%%%%%%%%%%%%%%%%%%%%%%%%%%%%%%%%%%%%%%%%%%%%%%%%%%%%%%%%%%%%%%%%%%%%%%%%%%%%%%%%%%%%%%%%%%%%%%%%%%%%%%%%%%%%%%%%%%%%%%%%%%%%%%%%%%%%%%%%%%%%%%%%%%%%%%%%%%%%%%%%%%%%%%%%%%%%%%%%%%%%%%%%%%%%%%%%%%%%%%%%%%%%%%%%%%%%%%
\section{KLR Algebras and Universal Specht Module}\label{sec:klr-specht}

In this section, we introduce KLR algebras, universal Specht modules, and the results we need to state our main result. For more details, readers are welcome to refer to \cite{kmr-universal-specht-type-A}.

\subsection{KLR algebras}
Let $\bk$ be a field, and fix $\alpha\in Q^+$ such that $\height(\alpha)=n$.

\begin{Definition}[\cite{khovanovlauda-klr-1,khovanovlauda-KLR-3,rouquier-2kacmoody}]\label{dfn-klr}
    The KLR algebra $R_{\alpha}$ of type $\Aone[e-1](e>2)$ is the unital $\bk$-algebra generated by the elements:
    \begin{equation}\label{klr-generating-elements}
        \{e(\bi)|\bi\in I^\alpha\}\cup \{y_1,\cdots,y_n\}\cup\{\psi_1,\cdots,\psi_{n-1}\}
    \end{equation}
    subject only to the following relations:
    \begin{align}
        e(\bi)e(\bj) &= \delta_{\bi,\bj} e(\bi), \qquad \sum\limits_{\bi \in I^\alpha} e(\bi) = 1 \label{eq:idempotents} \\ 
        %%%%%%%%%%%%%%%%%%%%%%%%%%
        y_r e(\bi) &= e(\bi) y_r,\qquad \psi_r e(\bi) = e(\sigma_r \bi) \psi_r \label{eq:commute_y_psi_e} \\
        %%%%%%%%%%%%%%%%%%%%%%%%%%%%%
        y_r y_s &= y_s y_r \label{eq:y_commute} \\
        %%%%%%%%%%%%%%%%%%%%%%%%%%%%%%%
        \psi_r y_s &= y_s \psi_r \quad \text{if } s \ne r, r+1 \label{eq:psi_y_commute} \\
        %%%%%%%%%%%%%%%%%%%%%%%%%%%%%%%
        \psi_r \psi_s &= \psi_s \psi_r \quad \text{if } |r - s| > 1 \label{eq:psi_commute_far} \\
        %%%%%%%%%%%%%%%%%%%%%%%%%%%%%%%
        \psi_r y_{r+1} e(\bi) &= (y_r \psi_r + \delta_{\bi_r, \bi_{r+1}}) e(\bi) \label{eq:psi_y_shift_left} \\
        %%%%%%%%%%%%%%%%%%%%%%%%%%%%%%%
        y_{r+1} \psi_r e(\bi) &= (\psi_r y_r + \delta_{\bi_r, \bi_{r+1}}) e(\bi) \label{eq:y_psi_shift_right} \\
        %%%%%%%%%%%%%%%%%%%%%%%%%%%%%%%
        \psi_r^2 e(\bi) &= Q_{\bi_r, \bi_{r+1}}(y_r, y_{r+1}) e(\bi) \label{eq:psi_square} \\
        %%%%%%%%%%%%%%%%%%%%%%%%%%%%%%%
        \psi_r \psi_{r+1} \psi_r e(\bi) &= \psi_{r+1} \psi_r \psi_{r+1} e(\bi) + Q_{\bi_r, \bi_{r+1}, \bi_{r+2}}(y_r, y_{r+1}, y_{r+2}) e(\bi) \label{eq:braid_relation}
    \end{align}
    where
    \[Q_{\bi_r, \bi_{r+1}}(y_r, y_{r+1})=\begin{cases}0 &\text{ if }\bi_r=\bi_{r+1}\\1& \text{ if }\bi_{r+1}\neq \bi_r,\bi_{r}\pm 1\\ y_{r+1}-y_r & \text{ if }\bi_r\rightarrow \bi_{r+1}\\ y_r-y_{r+1}&\text{ if }\bi_r\leftarrow \bi_{r+1}\end{cases}\]
    and 
    \[Q_{\bi_r, \bi_{r+1},\bi_{r+2}}(y_r, y_{r+1},y_{r+2})=\begin{cases}1 &\text{ if }\bi_r=\bi_{r+2}\to \bi_{r+1}\\-1& \text{ if }\bi_r=\bi_{r+2}\leftarrow \bi_{r+1}\\ 0&\text{ else }\end{cases}\]
\end{Definition}
Given any dominant weight $\Lambda\in P^+$, the corresponding cyclotomic KLR algebra $R^\Lambda_\alpha$ is generated by the same elements \autoref{klr-generating-elements} subject only to the above relations with the additional cyclotomic relations
\begin{equation}\label{eq:cyclotomic}
    y_1^{(\Lambda,\alpha_{\bi_1})}e(\bi)=0, \text{ where }\bi=(\bi_1,\cdots,\bi_n)
\end{equation}

Most importantly, $R_\alpha$ and $R^\Lambda_\alpha$ have $\mathbb{Z}$-gradings determined by setting $e(\bi)$ to be of degree $0$, $y_r$ of degree $2$, and $\psi_re(\bi)$ of degree $-a_{\bi_r,\bi_{r+1}}$ for all $r$ and $\bi\in I^\alpha$.

Each $w\in \Sym_n$ can be written as $w=\prod\limits_{1\leq j\leq m} \sigma_{i_j}$ for some $m$. If the expression is reduced (that is, of minimal length), then we define $\psi_w:=\psi_{i_1}\cdots\psi_{i_m}$. However, this definition depends on the choice of reduced expression of $w$, since \autoref{eq:braid_relation} breaks the usual braid relation. We therefore fix, for each $w\in \Sym_n$, a choice of reduced expression, and define $\psi_w$ with respect to that choice.

It is well-known \cite[Theorem 2.5]{khovanovlauda-klr-1} that $R_\alpha$ has a $\bk$-basis $\{\psi_w y_1^{a_1}\cdots y_n^{a_n}e(\bi)\}$ where $w\in\Sym_n$, $a_1,\cdots,a_n\in \mathbb{Z}_{\geq 0}$, $\bi\in I^\alpha$. 

Let $\kappa$ be a charge of $\Lambda$, and let $\blam\in\Par_{\alpha}$. Let $T$ be a row-standard $\blam$-tableau. Define $\psi^T:=\psi_{w^T}$, using our fixed choice of reduced expression for $w^T$.
%%%%%%%%%%%%%%%%%%%%%%%%%%%%%%%%%%%%%%%%%%%%%%%%%%%%%%%%%%%%%%%%%%%%%%%%%%%%%%%%%%%%%%%%%%%%%%%%%%%%%%%%%%%%%%%%%%%%%%%%%%%%%%%%%%%%%%%%%%%%%%%%%%%%%%%%%%%%%%%%%%%%%%%%%%%%%%%%%%%%%%%%%%%%%%%%%%%%%%%%%%%%%%%%%%%%%%%%%%%%%%%%%%%%%%%%%%%%%%%%%%%%%%%%%%%%%%%%%%%%%%%%%%%%%%%%%%%%%%%%%%%%%%%%%%%%%%%%%%%%%%%%%%%%%%%%%%%%%%%%%%%%%%%%%%%%%%%%%%%%%%%%%%%%%%%%%%%%%%%%%%%%%%%%%%%%%%%%%%%%%%%%%%%%%%%%%%%%%%%%%%%%%%%%%%%%%%%%%%%%%%%%%%%%%%%%%%%%%%%%%%%%%%%%%%%%%%%%%%%%%%%%%%%%%%%%%%%%%%%%%%%%%%%%%%%%%%%%%%%%%%%%%%%%%%%%%%%%%%%%%%%%%%%%%%%%%%%%%%%%%%%%%%%%%%%%%%%%%%%%%%%%%%%%%%%%%%%%%%%%%%%%%%%%%%%%%%%%%%%%%%%%%%%%%%%%%%%%%%%%%%%%%%%%%%%%%%%%%%%%%%%%%%%%%%%%%%%%%%%%%%
\subsection{Universal Specht module}\label{subsec:}
In this section, we fix $\alpha\in Q^+$, $\Lambda\in P^+$, and a charge $\kappa$ of $\Lambda$. We also fix $\blam \in \Par[\kappa]_\alpha$. Following \cite{kmr-universal-specht-type-A}, we define the Specht module $S^\blam$ over $R_\alpha$ in this section.

\begin{Remark}
    In \cite{kmr-universal-specht-type-A}, the quiver is oriented by $i\to i-1$, whereas our convention in \autoref{subsec:cartan-data} is the opposite. Nevertheless, we find the convention $i\to i+1$ more natural. Moreover, by replacing $i$ with $-i$ (modulo $e$), the arguments in \cite{kmr-universal-specht-type-A} that we use still apply. See also \cite[Lemma 3.2]{aip-rep-type-a}.
\end{Remark}

\begin{Definition}\label{dfn:garnir}
    A node $A=(m,r,c)\in[\blam]$ is a \textbf{Garnir node} of $\blam$ if $(m,r+1,c)\in[\blam]$. The \textbf{Garnir belt} of $A$ is the set $\mathcal{B}^A$ of nodes of $[\blam]$ consisting of $A$ and all nodes directly to the right of $A$, together with the node directly below $A$ and all nodes directly to the left of that node in the same component. Explicitly,
    \[
        \mathcal{B}^A=\{(m,r,z)\in[\blam]\mid c\le z\le \blam^{(m)}_r\}\cup \{(m,r+1,z)\in[\blam]\mid 1\le z\le c\}.
    \]
\end{Definition}
\begin{Definition}\label{dfn:garnir-tableau}
    Let $A\in[\blam]$ be a Garnir node. The \textbf{Garnir tableau} $G^A$ is the unique row-standard tableau satisfying:
    \begin{itemize}
        \item it agrees with $T^{\blam}$ on all nodes outside the Garnir belt $\mathcal{B}^A$,
        \item its entries in $\mathcal{B}^A$ increase from the bottom-left to the top-right.
    \end{itemize}
\end{Definition}
%%%%%%%%%%%%%%%%%
\begin{Example}\label{eg:garnir-tableau}
    As in \autoref{eg:tableau}, consider the partition $\lambda=(4,2)$. The Garnir nodes of $[\lambda]$ are $(1,1)$ and $A:=(1,2)$. The Garnir belt $\mathcal{B}^A$ consists of the following shaded nodes:
    \[
        \Tableau{~*~*~*~,*~*~}
    \]
    The Garnir tableau $G^A$ is the following tableau:
    \[
        \Tableau{1456,23}
    \]
    and $w^{G^A}=(\sigma_3\sigma_4\sigma_5)(\sigma_2\sigma_3\sigma_4)\in\Sym_6$.
\end{Example}
% To define the universal Specht module, for each Garnir node $A \in [\blam]$, we associate an element $g^A \in R_\alpha$, which yields the corresponding Garnir relation used in \autoref{dfn-universal-specht}. However, the definition of $g^A$ is lengthy, and in this paper we consider only the trivial case. Rather than presenting the full construction, we state the following result instead.

% \begin{Lemma}
% $g^A = \psi^{G^A}$ holds for any Garnir node $A$ of a hook partition $\lambda=(k,1^r)$. The same holds for any Garnir node of any partition in type~$A_\infty$ or type~$\Aone[e-1]$ when $e \gg 0$.
% \end{Lemma}
% \begin{proof}
% Let $\mathcal{D}^A$ be the set of minimal length left coset representatives of $\Sym_f\times \Sym_{k-f}$ defined in \cite[Section~5.2]{kmr-universal-specht-type-A}.  If $\mathcal{D}^A = \{1\}$, then $g^A = \psi^{G^A}$ follows directly from the setup in \cite[Definition~5.8]{kmr-universal-specht-type-A}. Observe that the Garnir nodes of a hook partition are precisely the nodes in the first column, excluding the bottom node. Since the second row of each Garnir belt in this case consists of exactly one node, we have $\mathcal{D} = \{1\}$ for every such Garnir node. The final statement holds because in type~$A_\infty$, a $\delta$-brick is infinitely long; hence, any finite partition cannot contain a $\delta$-brick, and thus any Garnir relation must be trivial. The same conclusion applies to $\Aone[e-1]$ when $e \gg 0$.
% \end{proof}
In this paper, we are focused on the following two cases:
\begin{equation}\label{concerned-case}
    \begin{cases}
        \text{(i) In type } \Aone[e-1], \text{ each component of } \blam \text{ is a hook partition,} \\
        \quad \text{i.e. } \blam^{(m)}=(k,1^r) \text{ for } k,r\in\mathbb{Z}_{\geq 1}; \\[1ex]
        \text{(ii) In type } A_{\infty}, \blam \text{ is an arbitrary } \ell\text{-partition.}
    \end{cases}
\end{equation}
We give the following specialized version of the definition of the universal Specht module from \cite[Definition~5.8]{kmr-universal-specht-type-A}.

\begin{Definition}\label{dfn-universal-specht}
    Fix $\alpha\in Q^+$ of height $n$, $\Lambda\in P^+$, and a charge $\kappa$ of $\Lambda$. Take $\blam\in\Par[\kappa]_{\alpha}$. Suppose we are in the situation of \autoref{concerned-case}. The \textbf{universal Specht module} $S^\blam$ is the cyclic $R_\alpha$-module generated by $z^\blam$, which is homogeneous of degree $\deg(T^\blam)$, subject only to the following relations:
    \begin{enumerate}
        \item $e(\bj)z^{\blam}= \delta_{\bj,\bi^{\blam}}z^{\blam}, \bj\in I^\alpha$;
        %%%%%%%%%%%
        \item $y_rz^{\blam}=0,  r=1,\dots,n$;
        %%%%%%%%%%%
        \item $\psi_r z^{\blam}=0$ whenever $r$ and $r+1$ appear in the same row of $T^{\blam}$;
        %%%%%%%%%%%
        \item $\psi^{G^A} z^{\blam}=0$ for all Garnir nodes $A \in [\blam]$. \label{garnir-relation}
    \end{enumerate}
\end{Definition}
%%%%%%%%%%
The last relation \autoref{garnir-relation} is called the \textbf{(trivial) Garnir relation} and $z^{\blam}$ is called the \textbf{standard cyclic generator} of $S^{\blam}$.

\begin{Remark}\label{rmk:trivial-garnir}
    \autoref{dfn-universal-specht} differs from \cite[Definition~5.8]{kmr-universal-specht-type-A} only in the fourth relation \autoref{garnir-relation}, where the latter uses the more general Garnir relation. However, since we restrict our attention to the cases in \autoref{concerned-case}, the two definitions agree because $\mathscr{D}^A = \{1\}$ for any Garnir node $A$ (cf.\ \cite[Definition~5.8]{kmr-universal-specht-type-A} and below).
\end{Remark}
In \autoref{dfn-universal-specht}, let $M^{\blam}$ be the cyclic $R_\alpha$-module with standard cyclic generator $z^{\blam}$ 
%\footnote{In \cite{kmr-universal-specht-type-A}, this cyclic generator of the permutation module is denoted by $m^{\blam}$. We abuse notation by continuing to write it as $z^{\blam}$.}
subject only to the first three relations. It is called the \textbf{permutation module}. Then $S^{\blam}$ is a quotient of $M^{\blam}$ by the Garnir relations, see \cite{kmr-universal-specht-type-A}. 

\medskip
In \cite[Section 3.6]{kmr-universal-specht-type-A}, a more detailed (and equivalent) construction of permutation modules using induction is provided. We describe this construction for partitions as follows:

Fix $\Lambda:=\Lambda_j$ a fundamental weight. For any one-row partition $\mu = (n)$ such that $\alpha_{\mu} = \alpha$, there exists a unique one-dimensional irreducible $R_\alpha$-module, denoted $L_\alpha$. It is generated by a vector $v$, defined by letting 
$e(\bi^{\mu})$ acts as the identity and all other generators act as zero.

Now consider an arbitrary partition (indeed, the construction works for any composition) $\lambda = (\lambda_1, \dots, \lambda_r) \in \Par_\alpha$. For each $1 \leq i \leq r$, define $\alpha(i) = \alpha_{\lambda_i}$ to be the positive root corresponding to the residue sequence of the $i$-th row of the partition $\lambda$ (with respect to $\Lambda$), and form the one-dimensional irreducible module $L_i := L_{\alpha(i)}$. We then define the module
\[
    L_1 \boxtimes L_2 \boxtimes \cdots \boxtimes L_r := L_1 \otimes_\bk L_2 \otimes_\bk \cdots \otimes_\bk L_r
\]
as a module over $R_{\alpha(1)} \otimes R_{\alpha(2)} \otimes \cdots \otimes R_{\alpha(r)}$. If we let $v_i$ be the cyclic generator of each $L_i$, then this module has cyclic generator $v_1 \boxtimes v_2 \boxtimes \cdots \boxtimes v_r$. Note that $\alpha = \alpha(1) + \cdots + \alpha(r)$.

By the standard theory of KLR algebras (see \cite{khovanovlauda-klr-1}), the algebra $R_{\alpha{(1)}} \otimes R_{\alpha{(2)}} \otimes \cdots \otimes R_{\alpha(r)}$ is a subalgebra of $R_\alpha$, and in fact, $R_\alpha$ is a free module over this subalgebra. Hence, we can apply the induction functor and define
\[
    M^\lambda := \Ind^{R_\alpha}_{R_{\alpha(1)} \otimes R_{\alpha(2)} \otimes \cdots \otimes R_{\alpha(r)}} L_1 \boxtimes L_2 \boxtimes \cdots \boxtimes L_r.
\]
The $R_\alpha$-module $M^\lambda$ is the permutation module corresponding to $\lambda$. It is straightforward to verify that $M^\lambda$ has the presentation stated below \autoref{dfn-universal-specht}; see \cite[Section 5.4]{kmr-universal-specht-type-A}. Clearly, this definition generalizes naturally to multipartitions.

\begin{Definition}\label{dfn:permutation-module-high-level}
Suppose $\blam = (\blam^{(1)} | \cdots | \blam^{(\ell)}) \in \Par_\alpha$ such that $\alpha(i) = \alpha_{\blam^{(i)}}$ for $1 \leq i \leq \ell$. Then we define
\[
    M^\blam := \Ind^{R_\alpha}_{R_{\alpha(1)} \otimes R_{\alpha(2)} \otimes \cdots \otimes R_{\alpha(\ell)}} M^{\blam^{(1)}} \boxtimes M^{\blam^{(2)}} \boxtimes \cdots \boxtimes M^{\blam^{(\ell)}}.
\]
\end{Definition}

For convenience, suppose that for each $1\le i\le \ell$ we have $\beta_i\in Q^+$ and an $R_{\beta_i}$-module $M_i$, and set $\beta:=\sum\limits_{1\le i\le \ell}\beta_i$. We then define
\begin{equation}\label{eq:def-of-circ}
    M_1\circ M_2\circ\cdots\circ M_{\ell}:=\Ind^{R_{\beta}}_{R_{\beta_1}\otimes R_{\beta_2}\otimes\cdots\otimes R_{\beta_{\ell}}}M_1\boxtimes M_2\boxtimes\cdots\boxtimes M_{\ell}.
\end{equation}
In particular, in \autoref{dfn:permutation-module-high-level},
$M^{\blam}=M^{\blam^{(1)}}\circ M^{\blam^{(2)}}\circ\cdots\circ M^{\blam^{(\ell)}}$.
%%%%%%%%%%%%%%%%%%%%%%%%%%%%%%%%%
\begin{Lemma}
    Let $\lambda \in \Par_\alpha$ be a partition of length $r$, and suppose $\Lambda = \Lambda_x$ is a fundamental weight. Let $\nu = (\lambda_1, \dots, \lambda_k)$ and $\mu = (\lambda_{k+1}, \dots, \lambda_r)$ for some $k \leq r$. Then
    \[
        M^\lambda \cong M^{(\nu | \mu)}, (\nu | \mu) \in \Par[\Lambda']_\alpha
    \]
    where $\Lambda'$ is determined by the charge $\kappa=(x, {x-k})$.
\end{Lemma}

\begin{proof}
This follows from the definition. One must carefully observe that the leading residue of each component is determined by the residue of $\lambda$. Specifically, the residue of the first node in the first row of the second component must equal that of the first node of the $(k+1)$-th row of $[\lambda]$.
\end{proof}
%%%%%%%%%%%%%%%%%%%%%%
\begin{Theorem}[{\cite[Theorem 5.6]{kmr-universal-specht-type-A}}]\label{thm:basis-permutation}
    The permutation module $M^{\blam}$ has a $\bk$-basis 
    \[
        \{\,\psi^T z^{\blam}\mid T\text{ is a row-standard }\blam\text{-tableau}\,\}
    \]
\end{Theorem}
\begin{Corollary}\label{cor:dim-permutation-module}
    For a partition $\lambda = (\lambda_1, \cdots, \lambda_r)$, let $n = \sum\limits_{1 \leq i \leq r} \lambda_i$. Then the dimension of $M^\lambda$ is given by:
    \[
        \dim M^\lambda = \frac{n!}{\lambda_1! \cdots \lambda_r!}.
    \]

\end{Corollary}
It turns out that the universal Specht modules of \autoref{dfn-universal-specht} are isomorphic to the graded Specht modules constructed in \cite{bkw-graded-specht} and the graded cell module constructed in \cite{humathas-graded-cellular}:
\begin{Theorem}[{\cite[Corollary 6.24]{kmr-universal-specht-type-A}}]\label{thm:basis-specht}
    There is a homogeneous degree $0$ isomorphism between $S^{\blam}$ and the graded cell module of
    $R^\Lambda_\alpha$ constructed in \textup{\cite{humathas-graded-cellular}}. Moreover,
    $S^{\blam}$ has a $\bk$-basis
    \[
        \{\psi^T z^{\blam} \mid T \in \Std(\blam)\}.
    \]
\end{Theorem}
One advantage of the construction of universal Specht modules is that the following result is immediate.
\begin{Theorem}[{\cite[Theorem 8.2]{kmr-universal-specht-type-A}}]\label{thm:kmr-high-level-specht}
    Suppose $\blam=(\blam^{(1)}\mid \cdots\mid \blam^{(\ell)})\in\Par_\alpha$, then 
    \[
        S^\blam\cong S^{\blam^{(1)}}\circ \cdots\circ S^{\blam^{(\ell)}}
    \]
    as (graded) $R_\alpha$-modules.
\end{Theorem}
%%%%%%%%%%%%%%%%%%%%%%%%%%%%%%%%%%%%%%%%%%%%%%%%%%%%%%%%%%%%%%%%%%%%%%%%%%%%%%%%%%%%%%%%%%%%%%%%%%%%%%%%%%%%%%%%%%%%%%%%%%%%%%%%%%%%%%%%%%%%%%%%%%%%%%%%%%%%%%%%%%%%%%%%%%%%%%%%%%%%%%%%%%%%%%%%%%%%%%%%%%%%%%%%%%%%%%%%%%%%%%%%%%%%%%%%%%%%%%%%%%%%%%%%%%%%%%%%%%%%%%%%%%%%%%%%%%%%%%%%%%%%%%%%%%%%%%%%%%%%%%%%%%%%%%%%%%%%%%%%%%%%%%%%%%%%%%%%%%%%%%%%%%%%%%%%%%%%%%%%%%%%%%%%%%%%%%%%%%%%%%%%%%%%%%%%%%%%%%%%%%%%%%%%%%%%%%%%%%%%%%%%%%%%%%%%%%%%%%%%%%%%%%%%%%%%%%%%%%%%%%%%%%%%%%%%%%%%%%%%%%%%%%%%%%%%%%%%%%%%%%%%%%%%%%%%%%%%%%%%%%%%%%%%%%%%%%%%%%%%%%%%%%%%%%%%%%%%%%%%%%%%%%%%%%%%%%%%%%%%%%%%%%%%%%%%%%%%%%%%%%%%%%%%%%%%%%%%%%%%%%%%%%%%%%%%%%%%%%%%%%%%%%%%%%%%%%%%%%%%%%
%%%%%%%%%%%%%%%%%%%%%%%%%%%%%%%%%%%%%%%%%%%%%%%%%%%%%%%%%%%%%%%%%%%%%%%%%%%%%%%%%%%%%%%%%%%%%%%%%%%%%%%%%%%%%%%%%%%%%%%%%%%%%%%%%%%%%%%%%%%%%%%%%%%%%%%%%%%%%%%%%%%%%%%%%%%%%%%%%%%%%%%%%%%%%%%%%%%%%%%%%%%%%%%%%%%%%%%%%%%%%%%%%%%%%%%%%%%%%%%%%%%%%%%%%%%%%%%%%%%%%%%%%%%%%%%%%%%%%%%%%%%%%%%%%%%%%%%%%%%%%%%%%%%%%%%%%%%%%%%%%%%%%%%%%%%%%%%%%%%%%%%%%%%%%%%%%%%%%%%%%%%%%%%%%%%%%%%%%%%%%%%%%%%%%%%%%%%%%%%%%%%%%%%%%%%%%%%%%%%%%%%%%%%%%%%%%%%%%%%%%%%%%%%%%%%%%%%%%%%%%%%%%%%%%%%%%%%%%%%%%%%%%%%%%%%%%%%%%%%%%%%%%%%%%%%%%%%%%%%%%%%%%%%%%%%%%%%%%%%%%%%%%%%%%%%%%%%%%%%%%%%%%%%%%%%%%%%%%%%%%%%%%%%%%%%%%%%%%%%%%%%%%%%%%%%%%%%%%%%%%%%%%%%%%%%%%%%%%%%%%%%%%%%%%%%%%%%%%%%%%%
\section{Specht Filtration in type \texorpdfstring{$\Aone[e-1]$}{A(1)\_\{e-1\}} for Hook Partitions}\label{sec:filtration}
%%%%%%%%%%%%%%%%%%%%%%%%%%%%%%%%%%%%%%%%%%%%%%%%%%%%%%%%%%%%%%%%%%%%%%%%%%%%%%%%%%%%%%%%%%%%%%%%%%%%%%%%%%%%%%%%%%%%%%%%%%%%%%%%%%%%%%%%%%%
In this section, we fix the quiver $\Aone[e-1]$ with $e > 2$, and take a positive root $\alpha \in Q^+$ and a dominant weight $\Lambda \in P^+$. Since our construction of the Specht filtration is component-by-component, see \autoref{thm:kmr-high-level-specht}, we may assume $\Lambda = \Lambda_x$ for some $0 \leq x \leq e-1$. 

Given a partition $\lambda = (\lambda_1, \dots, \lambda_r)\in\Par[\Lambda]_{\alpha}$, we define the corresponding permutation module $M^\lambda$. In general, there could be many Specht filtrations of $M^\lambda$. For example, let $\mu_i=(\lambda_i)\in\Par[\Lambda_{\res(i,1)}]$ for each $1\leq i\leq r$ and $\bmu=(\mu_1|\cdots|\mu_r)$, then (see \autoref{eq:def-of-circ})
\[
    M^\lambda\cong M^{\mu_1}\circ\cdots\circ M^{\mu_r}\cong S^{\mu_1}\circ \cdots\circ S^{\mu_r}\cong S^{\bmu}
\]
by \autoref{dfn:permutation-module-high-level} and \autoref{thm:kmr-high-level-specht}. However, this approach is not natural and proves to be unhelpful for our purposes. The reason is as follows: 

Although $M^\lambda$ is an $R_\alpha$-module, it is primarily used to construct the Specht module $S^\lambda$ over the cyclotomic algebra $R^\Lambda_\alpha$, where $\Lambda$ is fixed. Hence, a meaningful Specht filtration of $M^\lambda$ should be one with head $S^\lambda$. In other words, we aim to construct a Specht filtration for the submodule of $M^\lambda$ generated by the Garnir relations associated to $[\lambda]$.

In general, constructing such a filtration is a difficult problem. However, for hook partitions, the situation simplifies considerably. One of the key reasons is that in this case, all Garnir relations are trivial as mentioned in \autoref{rmk:trivial-garnir}.
%%%%%%%%%%%%%%%%%%%%%%%%%%%%%%%%%%%%%%%%%%%%%%%%%%%%%%%%%%%%%%%%%

%%%%%%%%%%%%%%%%%%%%%%%%%%
\medskip
Throughout this section, we fix a hook partition $\lambda = (k, 1^r)$. For each $1 \leq i \leq r$, let $A_i = (i, 1) \in [\lambda]$. One checks easily that $\{A_i \mid 1 \leq i \leq r\}$ is the set of all Garnir nodes in $[\lambda]$. We set
\[
    \psi^{A_i} := \psi^{G^{A_i}} \quad \text{for all } 1 \le i \le r.
\]

\begin{Lemma}\label{lm:garnir-relation-hook-partition}
    Let $v=z^{\lambda}$ be the standard cyclic generator of the permutation module $M^\lambda$ over $R_\alpha$, where $\alpha = \alpha_{\lambda}$. Then:
    \[
    \psi^{A_i}v =
    \begin{cases}
        \psi_1 \psi_2 \cdots \psi_k e(\bi^{\lambda})v & \text{if } i = 1, \\
        \psi_{k+i-1} e(\bi^{\lambda})v & \text{if }2\leq i\leq r
    \end{cases}
    \]
\end{Lemma}

\begin{Theorem}[\textbf{Specht Filtration}]\label{thm-specht-filtration-hook}
    Fix the quiver $\Aone[e-1]$ ($e>2$) or $A_\infty$, take $\alpha\in Q^+$, $\Lambda=\Lambda_x$ and a hook partition $\lambda=(k,1^r)$ such that $\alpha_{\lambda}=\alpha$. 
    A Specht filtration of the permutation module $M^\lambda$ is given by the following chain of $R_\alpha$-modules:
    \[
    M^\lambda = M_0 \supsetneq M_1 \supsetneq M_2 \supsetneq \cdots \supsetneq M_r \supsetneq M_{r+1} = 0
    \]
    where, for $1 \leq i \leq r$, the module $M_i$ is the submodule of $M^{\lambda}$ generated by $\{\psi^{A_i}v, \dots, \psi^{A_r}v\}$, with $v=z^{\lambda}$ being the standard cyclic generator of $M^{\lambda}$ over $R_{\alpha}$. Moreover,
    \[
    M_i / M_{i+1} \cong S^{\lambda_i}
    \]
    where
    \[
    \lambda_i =
    \begin{cases}
        \lambda = (k, 1^r) & \text{if } i = 0, \\
        (k+1, 1^{r-1}) & \text{if } i = 1, \\
         \big(k | \underbrace{ 1 | 1 | \cdots | 1}_{i-2 \text{ times}} 
 | (2, 1^{r-i})\big)& \text{if } 2 \leq i \leq r, 
        %(k) | \underbrace{ (1) | (1) | \cdots | (1)}_{r-2 \text{ times}} | (2) & \text{if } i = r.
    \end{cases}
    \]
    and $S^{\lambda_i}$ is the Specht module over $R^{\Lambda(i)}_\alpha$ where $\Lambda(i)$ is determined by the charge:
    \[
    \kappa(i) =
    \begin{cases}
        x & \text{if } i = 0, \\
        {x-1} & \text{if } i = 1, \\
         ({x}, {x-1},\cdots,{x-i+2},{x-i})& \text{if } 2 \leq i \leq r
    \end{cases}
    \]
\end{Theorem}

%%%%%%%%%%%%%%%%%%%%%%%%%%
Before giving the proof, we show the construction in an example.
\begin{Example}
     Take quiver $\Aone[9]$ with vertex set $\{0,1,2,\cdots,9\}$ and $\Lambda=\Lambda_0$, consider \[\beta=\alpha_0+\alpha_1+\alpha_2+\alpha_3+\alpha_5+\alpha_6+\alpha_7+\alpha_8+\alpha_9\] and $\lambda=(4,1^5)$.  
    
    We know $\dim M^\lambda=15120$. The partitions that appear in the Specht filtration of  $M^\lambda$ are listed below, with their corresponding Young diagrams filled with residues.
    \begin{enumerate}
        \item 
        $\lambda_0 = \lambda = (4,1^5)$
        \[
            \Tableau{0123,9,8,7,6,5} \qquad 
            \dim S^{\lambda_0} = 56 \qquad 
            \Lambda(0) = \Lambda_0
        \]

        \item 
        $\lambda_1 = (5,1^4)$
        \[
            \Tableau{90123,8,7,6,5} \qquad 
            \dim S^{\lambda_1} = 70 \qquad 
            \Lambda(1) = \Lambda_9
        \]

        \item 
        $\lambda_2 = (4|2,1^3)$
        \[
            \Multitableau{0123|89,7,6,5} \qquad 
            \dim S^{\lambda_2} = 504 \qquad 
            \Lambda(2) = \Lambda_0 + \Lambda_8
        \]

        \item 
        $\lambda_3 = (4|1|2,1^2)$
        \[
            \Multitableau{0123|{9}|78,6,5} \qquad 
            \dim S^{\lambda_3} = 1890 \qquad 
            \Lambda(3) = \Lambda_0 + \Lambda_9 + \Lambda_7
        \]

        \item 
        $\lambda_4 = (4|1|1|2,1)$
        \[
            \Multitableau{0123|9|8|67,5} \qquad 
            \dim S^{\lambda_4} = 5040 \qquad 
            \Lambda(4) = \Lambda_0 + \Lambda_9 + \Lambda_8 + \Lambda_6
        \]
        
        \item 
        $\lambda_5=(4|1|1|1|2)$
        \[
            \Multitableau{0123|9|8|7|56},\quad 
            \dim S^{\lambda_5} = 7560,\quad 
            \Lambda(5) = \Lambda_0 + \Lambda_9 + \Lambda_8 + \Lambda_7 + \Lambda_5
        \]
    \end{enumerate}
    It is easy to verify that:
    $\sum\limits_{0\leq i\leq 5}\dim S^{\lambda_i}=56+70+504+1890+5040+7560=15120=\dim M^\lambda$.
\end{Example}
%%%%%%%%%%%%%%%%%%%%%%%%%%%%%%%%%%%%%%%%%%%%%%%%%%%%%%%%%%%%%%%%%%%%%%%%%%%%%%%%%%%%%%%%%%%%%%%%%%%%%%%%%%%%%%%%%%%%%%%%%%%%%%%%%%
The equality of dimensions plays a crucial role in our proof, so we begin by establishing this fact.
\begin{Lemma}\label{lm:equality-dimension-specht-permutation}
    In type $\Aone[e-1]$ or type $\text{A}_\infty$, suppose $\lambda=(k,1^r)$ is a hook partition, then $\dim M^\lambda=\sum\limits_{i=0}^{r}\dim S^{\lambda_i}$.
\end{Lemma}
\begin{proof}
    We use the \emph{hook length formula} to compute the dimension of Specht modules; see, for instance, \cite[Theorem 4.2.14]{csst-rep-symmetric}. Then
    \[\dim M^\lambda=\frac{(k+r)!}{k!},\] 
    \[\dim S^{\lambda_0}=\dim S^{\lambda}=\frac{(k+r)!}{(k+r)r!(k-1)!}=\frac{(k+r)!k}{(k+r)r!k!},\] 
    \[\dim S^{\lambda_1}=\frac{(k+r)!}{(k+r)k!(r-1)!}=\frac{(k+r)!r}{(k+r)k!r!}.\]
    %%%%%%%%%%%%%%%%%%%%%%%%%%%%%%%%%%
    For $2\leq i\leq r-1$, it is easy to see:
    \begin{align*}
    \dim S^{\lambda_i}
    &= \binom{k+r}{k}\binom{r}{1}\binom{r-1}{1}\cdots\binom{r-i+3}{1}(r-i+1)\\
    %&=\frac{(k+r)!}{k!r!}\frac{r(r-1)\cdots(r-i+3)(r-i+1)}{1}\\
    &=\frac{(k+r)!}{k!r!}\frac{r(r-1)\cdots(r-i+3)(r-i+2)(r-i+1)}{r-i+2}\\
    &=\frac{(k+r)!}{k!r!}\frac{r!}{(r-i)!(r-i+2)}\\
    &=\frac{(k+r)!}{k!(r-i)!(r-i+2)}\\
    &=\frac{(k+r)!(r-i+1)}{k!(r-i+2)!}.
    \end{align*}
    For $S^{\lambda_r}$, we know:
    \[\dim S^{\lambda_r}=\frac{(k+r)!}{k!2}.\]
    Hence we have the following:
    \begin{align*}
            \sum\limits_{0\leq i\leq r}\dim S^{\lambda_i}&=\frac{(k+r)!}{k!}\left( \frac{k}{(k+r)r!}+\frac{r}{(k+r)r!}+\sum\limits_{2\leq i\leq r-1}\frac{r-i+1}{(r-i+2)!}+\frac{1}{2} \right)\\
            %&=\frac{(k+r)!}{k!}\left( \frac{1}{r!}+\sum\limits_{2\leq i\leq r-1}\frac{r-i}{(r-i+2)!}+\sum\limits_{2\leq i\leq r-1}\frac{1}{(r-i+2)!}+\frac{1}{2} \right)\\
            &=\frac{(k+r)!}{k!}\left( \frac{1}{r!}+\sum\limits_{2\leq i\leq r-1}\frac{r-i+2}{(r-i+2)!}-\sum\limits_{2\leq i\leq r-1}\frac{1}{(r-i+2)!}+\frac{1}{2} \right)\\
            &=\frac{(k+r)!}{k!}\left( \frac{1}{r!}+\sum\limits_{2\leq t\leq r-1}\frac{1}{t!}-\sum\limits_{3\leq t\leq r}\frac{1}{t!}+\frac{1}{2} \right)\\
            &=\frac{(k+r)!}{k!}=\dim M^\lambda.
    \end{align*}
\end{proof}
\autoref{lm:equality-dimension-specht-permutation} shows that the filtration of \autoref{thm-specht-filtration-hook} is well-behaved at the level of vector spaces.
%%%%%%%%%%%%%%%%%%%%%%%%%%%%%%%%%%%%%%%%%%%%%%%%
From this point onward, we fix the notation as in \autoref{thm-specht-filtration-hook} and assume $\lambda = (k,1^r)$ with $k>1$ and $r>1$ to avoid trivial cases. The case $r=1$ is analogous, and can be handled by the same argument as in \autoref{lm-proof-i=1}.
\begin{Lemma}\label{lm-relation-ypsi-psipsi}
    Let $v=z^\lambda$ be the standard cyclic generator of $M^\lambda$ and $\bi=\bi^\lambda=\res(T^\lambda)$, we have:
    \begin{align}
        y_j\left(\psi_1\psi_2\cdots\psi_k e(\bi) v\right) &= 0, && \text{for } 1 \leq j \leq k+r, \tag{a} \label{eq:y-i=1} \\
        %%%%%%%%%%%%%%%
        y_j\left(\psi_{k+i-1} e(\bi) v\right) &= 0, && \text{for } 2 \leq i \leq r,\; 1 \leq j \leq k+r, \tag{b} \label{eq:y-i>1} \\
        %%%%%%%%%%%%%%%
        \psi_j\left(\psi_1\psi_2\cdots\psi_k e(\bi) v\right) &= 0, && \text{for } 1 \leq j \leq k, \tag{c} \label{eq:psi-i=1} \\
        %%%%%%%%%%%%%%%
        \psi_j\left(\psi_{k+i-1} e(\bi) v\right) &= 0, && \text{for } 2 \leq i \leq r,\; 1 \leq j \leq k-1. \tag{d} \label{eq:psi-i>1}
    \end{align}
\end{Lemma}
\begin{proof}
    To show \autoref{eq:y-i=1}, if $j=1$, then the tableau of $\sigma_2\sigma_3\cdots\sigma_k T^\lambda$ is the following:
    \begin{center}
        \Tableau[xscale=1.4,dotted cols={4,5,6},dotted rows={4,5,6}]{134000k{k{+}1},2,{k{+}2},0,0,0,{k{+}r}}
    \end{center}
    Hence we have:
    \[y_1\psi_1\bigl(\psi_2\cdots\psi_k e(\bi) v\bigr)=\psi_1y_2\bigl(\psi_2\cdots\psi_ke(\bi)v\bigr).\]
    For each $2\leq i\leq k-1$, the tableau $\sigma_{i+1}\sigma_{i+2}\cdots\sigma_{k}T^\lambda$ is the following:
    \begin{center}
        \Tableau[xscale=1.4,dotted cols={3,4,7,8},dotted rows={4,5}]{1200{i}{i{+}2}00k{k{+}1},{i{+}1},{k{+}2},0,0,{k{+}r}}
    \end{center}
    Let $\delta_i=\delta_{\res(1,i),\res(2,1)}\in\{0,1\}$, then we have:
    \begin{align*}
        \psi_1\cdots\psi_{i-1}\bigl(y_i\psi_{i}\bigr)\psi_{i+1}\cdots\psi_{k}v
        &=\psi_1\cdots\psi_{i-1}\bigl(\psi_{i}y_{i+1}-\delta_i\bigr)\psi_{i+1}\cdots\psi_{k}e(\bi)v\\
        &=\psi_1\cdots\psi_{i-1}\psi_{i}y_{i+1}\psi_{i+1}\cdots\psi_{k}e(\bi)v
        -\delta_i\psi_1\cdots\psi_{i-1}\psi_{i+1}\cdots\psi_{k}e(\bi)v\\
        &=\psi_1\cdots\psi_{i-1}\psi_{i}y_{i+1}\psi_{i+1}\cdots\psi_{k}e(\bi)v
        -\delta_i\psi_1\cdots\psi_{i-2}\psi_{i+1}\cdots\psi_{k}\psi_{i-1}e(\bi)v\\
        &=\psi_1\cdots\psi_{i-1}\psi_{i}y_{i+1}\psi_{i+1}\cdots\psi_{k}e(\bi)v.
    \end{align*}
    Hence we have:
    \begin{align*}
        y_1\psi_1\psi_2\cdots\psi_k e(\bi) v&=\psi_1y_2\psi_2\cdots\psi_ke(\bi)v\\
        &=\psi_1\psi_2\cdots y_i\psi_i\cdots\psi_ke(\bi)v\\
        &=\psi_1\psi_2\cdots y_k\psi_k e(\bi)v\\
        &=\psi_1\psi_2\cdots\psi_{k-1}\bigl(\psi_k y_{k+1}-\delta_k\bigr) e(\bi)v\\
        &=\psi_1\psi_2\cdots\psi_{k-1}\psi_k y_{k+1} e(\bi)v=0.
    \end{align*}
    If $1<j\le k+1$, we have:
    \begin{align*}
        y_j\psi_1\psi_2\cdots\psi_k e(\bi) v&=\psi_1\psi_2\cdots \psi_{j-2}y_j\psi_{j-1}\cdots\psi_ke(\bi)v\\
        &=\psi_1\psi_2\cdots \psi_{j-2}(\psi_{j-1}y_{j-1}+\delta_{j-1})\psi_{j}\cdots\psi_k e(\bi)v\\
        &=\psi_1\psi_2\cdots \psi_{j-2}\psi_{j-1}y_{j-1}\psi_{j}\cdots\psi_k e(\bi)v\\
        &=\psi_1\psi_2\cdots \psi_{j-2}\psi_{j-1}\psi_{j}\cdots\psi_ky_{j-1}e(\bi)v=0.
    \end{align*}
    If $j>k+1$, then $y_j$ commutes with $\psi_1\psi_2\cdots \psi_k$, and hence \autoref{eq:y-i=1} is trivial.
    
%%%%%%%%%%%%%%%%%%%%%%%%%%%%%%%%%%%%%%%%%%%%%
    \medskip
    To show \autoref{eq:y-i>1}, if $j\neq k+i-1,k+i$, then $y_j$ commutes with $\psi_{k+i-1}$ and it follows. As $\bi_{k+i-1}=\res(i,1)\leftarrow\res(i+1,1)=\bi_{k+i}$, we have:
    \[y_{k+i-1}\psi_{k+i-1}e(\bi)v=\psi_{k+i-1}y_{k+i}e(\bi)v=0,\]
    \[y_{k+i}\psi_{k+i-1}e(\bi)v=\psi_{k+i-1}y_{k+i-1}e(\bi)v=0.\]
    %%%%%%%%%%%%%%%%%%%%%%%%%%%%%%%%%%%%%%%%%%%%%%%%
    \medskip
    To show \autoref{eq:psi-i=1}, if $j=1$,
    \begin{align*}
        \psi_1^2\psi_2\cdots\psi_k e(\bi) v 
        &=(y_1-y_2)\psi_2\cdots\psi_k e(\bi) v\\
        &=-y_2\psi_2\cdots\psi_k e(\bi) v+\psi_2\cdots\psi_k y_1e(\bi) v\\
        &=-y_2\psi_2\cdots\psi_k e(\bi) v\\
        &=-\psi_2y_3\psi_3\cdots \psi_ke(\bi)v\\&=0.
    \end{align*}
    The last equality holds by above argument in \autoref{eq:y-i=1}. The second last equality holds because $\delta_{2}=\delta_{\res(1,2),\res(2,1)}=0$ since $e>2$.% This implies in type C, it is a big issue.

    If $j=2$, we use \autoref{eq:braid_relation} and that $\delta_1=\delta_{\res(1,1),\res(2,1)}=0$:
    \begin{align*}
        \psi_2\psi_1\psi_2\psi_3\cdots\psi_k e(\bi)v
        &=(\psi_1\psi_2\psi_1-\delta_{1})\psi_3\cdots\psi_k e(\bi)v\\
        &=\psi_1\psi_2\psi_1\psi_3\cdots\psi_k e(\bi)v\\
        &=\psi_1\psi_2\psi_3\cdots\psi_k\psi_1 e(\bi)v\\
        &=0.
    \end{align*}
    
    If $2<j\leq k$, then we use \autoref{eq:braid_relation}:
    \begin{align*}
        \psi_j\psi_1\cdots\psi_ke(\bi)v
        &=\psi_1\cdots\psi_{j}\psi_{j-1}\psi_{j}\cdots\psi_k e(\bi)v\\
        &=\psi_1\cdots(\psi_{j-1}\psi_{j}\psi_{j-1}-\delta_{j-1})\cdots\psi_k e(\bi)v\\
        &=\psi_1\cdots\psi_{j-1}\psi_{j}\psi_{j-1}\cdots\psi_k e(\bi)v-\delta_{j-1}\psi_1\cdots\psi_{j-2}\psi_{j+1}\cdots\psi_k e(\bi)v\\
        &=\psi_1\cdots\psi_{j-1}\psi_{j}\cdots\psi_k \psi_{j-1}e(\bi)v-\delta_{j-1}\psi_1\cdots\psi_{j+1}\cdots\psi_k\psi_{j-2} e(\bi)v\\
        &=0.
    \end{align*}
    \medskip
    To show \autoref{eq:psi-i>1}, notice that $\psi_j$ commutes with $\psi_{k+i-1}$ and $\psi_je(\bi)v=0$ for any $2 \leq i \leq r, 1 \leq j \leq k-1$. 
\end{proof}
%%%%%%%%%%%%%%%%%%%%%%%%%%%%%%%%%%%%%%%%%%%
\begin{Lemma}\label{lm-proof-i=1}
    There is a surjective $R_\alpha$-homomorphism from $S^{\lambda_1}$ to $M_1/M_2$, which maps the standard cyclic generator $w$ of $S^{\lambda_1}$ to $\psi^{A_1}v+M_2$.
\end{Lemma}

\begin{proof}
    Let $w$ be the standard cyclic generator of the Specht module $S^{\lambda_1}$. It has the following presentation:
    \begin{align}
        e(\bj')w &= \delta_{\bj,\bj'}e(\bj)w, && \text{where } \bj = \res(T^{\lambda_1}) \tag{1} \label{eq:1} \\[1ex]
        %%%%%%%%%%%%%%%%%%
        y_j e(\bj) w &= 0, && \text{for } 1 \leq j \leq k + r \tag{2} \label{eq:2} \\[1ex]
        %%%%%%%%%%%%%%%%%%
        \psi_j e(\bj) w &= 0, && \text{for } 1 \leq j \leq k \tag{3} \label{eq:3} \\[1ex]
        %%%%%%%%%%%%%%%%%%
        \psi_1 \psi_2 \cdots \psi_{k+1} e(\bj) w &= 0 \tag{4} \label{eq:4} \\[1ex]
        %%%%%%%%%%%%%%%%%%%
        \psi_{k+i} e(\bj) w &= 0, && \text{for } 2 \leq i \leq r - 1 \tag{5} \label{eq:5}
    \end{align}

    We verify that $M_1/M_2$ with the standard cyclic generator $\psi^{A_1}v+M_2=\psi_1\psi_2\cdots\psi_k e(\bi^{\lambda}) v+M_2$ satisfies these relations, hence there is a  surjective $R_\alpha$-homomorphism $\phi_1$ from $S^{\lambda_1}$ to $M_1/M_2$ mapping $w$ to $\psi^{A_1}v+M_2$.

    %Let the Garnir nodes of $[\lambda_1]$ be $B_1,\cdots,B_{r-1}$ where $B_i:=(i,1)\in [\lambda_1]$. 
    
    The \autoref{eq:1} relation is clear since $\res_{\Lambda_x}(G^{A_1})=\res_{\Lambda_x}(\sigma_1\sigma_2\cdots\sigma_k T^\lambda)=\res_{\Lambda_{x-1}}(T^{\lambda_1})$.%need to define the Garnir tableau as well as the notation $\res[\Lambda]$

    The \autoref{eq:2} and \autoref{eq:3} relations hold by \autoref{eq:y-i=1} and \autoref{eq:psi-i=1} from \autoref{lm-relation-ypsi-psipsi}.

    For \autoref{eq:4}, notice that $T_j:=\sigma_j\sigma_{j+1}\cdots \sigma_{k+1}(\sigma_{j-1}\sigma_j\cdots\sigma_kT^{\lambda})$ $(j\geq 3)$ is the following tableau:
    \begin{center}
        \Tableau[xscale=1.4,dotted cols={3,4,7,8},dotted rows={6,7}]{1200{j{-}2}{j{+}1}00{k{+}2},{j{-}1},j,{k{+}3},{k{+}4},0,0,{k{+}r}}
    \end{center}
    Let $\bj':=\res(T_j)$. 
    It is easy to see $\res_{j-2}(T_j)= \res_{j}(T_j)\rightarrow \res_{j-1}(T_j)$ if and only if $\res(1,j-2)=\res(3,1)$ and this is the only case for which $Q_{\bj_{j-2}',\bj_{j-1}',\bj_{j}'}(y_{j-2},y_{j-1},y_{j})\neq 0$. Let $\delta^j=\delta_{\res(1,j-2),\res(3,1)}$, then 
    we have:
    \[ 
        Q_{\bj_{j-2}',\bj_{j-1}',\bj_{j}'}(y_{j-2},y_{j-1},y_{j})=\delta^j.
    \]
    Notice that $\delta^3=0$. We keep applying the relations \autoref{eq:braid_relation}:
    \begin{align*}
        \psi_1 \psi_2 \psi_3 \psi_4 \cdots \psi_{k+1} \bigl(\psi_1 \psi_2 \cdots \psi_k v + M_2\bigr)
        &= (\psi_1 \psi_2 \psi_1) \psi_3 \psi_4 \cdots \psi_{k+1} \bigl(\psi_2 \psi_3 \cdots \psi_k v + M_2\bigr) \\
        &= (\psi_2 \psi_1 \psi_2) \psi_3 \psi_4 \cdots \psi_{k+1} \bigl(\psi_2 \psi_3 \cdots \psi_k v + M_2\bigr)\\
        &= \psi_2\psi_1(\psi_2\psi_3\psi_2)\psi_4\cdots\psi_{k+1}\bigl(\psi_3\cdots\psi_kv+M_2\bigr)\\
        &= (\psi_2\psi_1)(\psi_3\psi_2\psi_3+\delta^4)\psi_4\cdots\psi_{k+1}\bigl(\psi_3\cdots\psi_kv+M_2\bigr)\\
        &= (\psi_2\psi_1)(\psi_3\psi_2\psi_3)\psi_4\cdots\psi_{k+1}\bigl(\psi_3\cdots\psi_kv+M_2\bigr)\\
        &=(\psi_2\psi_1)(\psi_3\psi_2)\cdots(\psi_{j-1}\psi_{j}\psi_{j-1})\psi_{j+1}\cdots\psi_{k+1}\bigl(\psi_{j}\cdots\psi_ke(\bi)v+M_2\bigr)\\
        &=(\psi_2\psi_1)(\psi_3\psi_2)\cdots(\psi_{j}\psi_{j-1}\psi_{j}+\delta^{j+1})\psi_{j+1}\cdots\psi_{k+1}\bigl(\psi_{j}\cdots\psi_ke(\bi)v+M_2\bigr)\\
        &=(\psi_2\psi_1)(\psi_3\psi_2)\cdots(\psi_{j}\psi_{j-1}\psi_{j})\psi_{j+1}\cdots\psi_{k+1}\bigl(\psi_{j}\cdots\psi_ke(\bi)v+M_2\bigr)\\
        &=\cdots\\
        &=(\psi_2\psi_1)(\psi_3\psi_2)\cdots(\psi_{k-1}\psi_{k-2})(\psi_k\psi_{k-1})\psi_{k}\psi_{k+1}\psi_{k}e(\bi)v+M_2\\
        &=(\psi_2\psi_1)(\psi_3\psi_2)\cdots(\psi_{k-1}\psi_{k-2})(\psi_k\psi_{k-1})(\psi_{k+1}\psi_k\psi_{k+1}+\delta^{k+2}) e(\bi)v+M_2\\
        &=(\psi_2\psi_1)(\psi_3\psi_2)\cdots(\psi_{k-1}\psi_{k-2})(\psi_k\psi_{k-1})\psi_{k+1}\psi_k\psi_{k+1}e(\bi)v+M_2\\
        &=0+M_2.
    \end{align*}
    All the expressions with $\delta^j$ $(4\leq j\leq k+2)$ vanish because the earlier term $\psi_{j-3}$ commutes with the terms to its right and kills $v$.
    The last equality holds because the Garnir relation $\psi^{A_2}v=\psi_{k+1}e(\bi^{\lambda})v\in M_2$.

    For \autoref{eq:5}, as already noted in \autoref{lm:garnir-relation-hook-partition}, we have
    \[
        \psi_{k+i}\psi^{A_1}v=\psi^{A_1}\psi_{k+i}v=\psi^{A_1}\psi^{A_{i+1}}v\in M_2,
    \]
    where the first equality holds because $\psi_{k+i}$ commutes with $\psi_1\cdots\psi_k$ for $2\le i\le r-1$.
\end{proof}
%%%%%%%%%%%%%%%%%%%%%%%%%%%%%%%%%%%%%%%%%%%%%%%%
\begin{Lemma}\label{lm-proof-r>i>1}
    For each $2\leq i\leq r-1$ there is a canonical surjective $R_\alpha$-homomorphism from $S^{\lambda_i}$ to $M_i/M_{i+1}$, which maps the standard cyclic generator $w$ of $S^{\lambda_i}$ to $\psi^{A_i}v+M_{i+1}$.
\end{Lemma}
\begin{proof}
    Let $w$ be the standard cyclic generator of the Specht module $S^{\lambda_i}$ and let $\bi:=\bi^{\lambda}$. Then it has the following presentation:
    \begin{align}
        e(\bi')w &= \delta_{\bj,\bi'}w, && \text{where } \bj = \res(T^{\lambda_i}) \tag{i} \label{eq:i} \\[1ex]
        %%%%%%%%%%%%%%%%%%%%
        y_j e(\bj) w &= 0, && \text{for } 1 \leq j \leq k + r \tag{ii} \label{eq:ii} \\[1ex]
        %%%%%%%%%%%%%%%%%%%
        \psi_j e(\bj) w &= 0, && \text{for } 1 \leq j \leq k-1 \text{ or } j = k+i-1 \tag{iii} \label{eq:iii} \\[1ex]
        %%%%%%%%%%%%%%%%%%%%%
        \psi_{k+i-1}\psi_{k+i} e(\bj) w &= 0 \tag{iv} \label{eq:iv} \\[1ex]
        %%%%%%%%%%%%%%%%%%%%%
        \psi_{k+i+j} e(\bj) w &= 0, && \text{for } 1 \leq j \leq r-1-i \tag{v} \label{eq:v}
    \end{align}

    We verify that $M_{i}/M_{i+1}$ with the standard cyclic generator $\psi^{A_i}v+M_{i+1}=\psi_{k+i-1}e(\bi) v+M_{i+1}$ satisfies these relations, hence there is a surjective $R_\alpha$- homomorphism $\phi_i$ from $S^{\lambda_i}$ to $M_i/M_{i+1}$ mapping $w$ to $\psi^{A_i}v+M_{i+1}$.

    %Let the Garnir nodes of $[\lambda_1]$ be $B_1,\cdots,B_{r-1}$ where $B_i:=(i,1)\in [\lambda_1]$. 
    \autoref{eq:i} can be shown by the same argument as in \autoref{lm-proof-i=1}.

    \autoref{eq:ii} and \autoref{eq:iii} are satisfied by \autoref{eq:y-i>1} and \autoref{eq:psi-i>1} from \autoref{lm-relation-ypsi-psipsi}.

    \autoref{eq:iv} can be verified by applying the relation \autoref{eq:braid_relation}:
    \[\psi_{k+i-1}\psi_{k+i}e(\bj)\bigl(\psi_{k+i-1}e(\bi) v+M_{i+1}\bigr)=\psi_{k+i}\psi_{k+i-1}\psi_{k+i}e(\bi) v+M_{i+1}=0+M_{i+1}.\]
    The first equality holds because $\bi_{k+i-1}\neq \bi_{k+i+1}$ since $e>2$.

    For the remaining relations from \autoref{eq:v}, use the fact $\psi_{k+i-1}$ commutes with $\psi_{k+i+j}$ for $1\leq j\leq r-1-i$:
    \[\psi_{k+i+j}e(\bj)\bigl(\psi_{k+i-1}v+M_{i+1}\bigr)=\psi_{k+i-1}\psi_{k+i+j}v+M_{i+1}.\]
    As $R_{\alpha}\{\psi_{k+i+1}v,\cdots,\psi_{k+r-1}v\}= M_{i+1}$, the conclusion follows.
\end{proof}
%%%%%%%%%%%%%%%%%%%%%%%%%%%%%%%%%%%%%%%%%%%%%%%%
\begin{Lemma}\label{lm-proof-i=r}
    There is a canonical surjective $R_\alpha$-homomorphism from $S^{\lambda_r}$ to $M_r$, which maps the standard cyclic generator $w$ of $S^{\lambda_r}$ to $\psi^{A_r}v$.
\end{Lemma}
\begin{proof}
    The proof is almost the same as \autoref{lm-proof-r>i>1}, but easier because there is no Garnir relation this time.
\end{proof}
%%%%%%%%%%%%%%%%%%%%%%%%%%%%%%%%%%%%%%%%%%%%%%%%
Now we can prove the main theorem.
 \begin{proof}[Proof of \autoref{thm-specht-filtration-hook}]
    The $i=0$ case is trivial since $M_1$ is generated by the Garnir relations and we know that
    \[
    M^\lambda/M_1 \cong S^\lambda
    \]
    as $R_\alpha$-modules.

    For $1 \leq i \leq r$, we have constructed a surjective $R_\alpha$-homomorphism from $S^{\lambda_i}$ to $M_i/M_{i+1}$ by \autoref{lm-proof-i=1}, \autoref{lm-proof-r>i>1} and \autoref{lm-proof-i=r}. By \autoref{lm:equality-dimension-specht-permutation}, since the sum of the dimensions of $S^{\lambda_i}$ for $0 \leq i \leq r$ equals the dimension of $M^\lambda$, all these surjective homomorphisms must be isomorphisms. Hence, the proof is complete.
\end{proof}
\begin{Remark}
See \cite{bkw-graded-specht} or \cite{kmr-universal-specht-type-A} for the definition of degree of tableaux. One might expect that the Specht filtration satisfies the property that the isomorphism between $S^{\lambda_i}$ and $M_i / M_{i+1}$ is homogeneous of degree zero. However, this is not true in general. Indeed, it is easy to see that, as graded $R_\alpha$-modules, we have:
$S^{\lambda_i}\<-\deg T^{\lambda_i}+\deg T^{\lambda}+\deg \psi^{A_i}e(\bi)\>\cong M_i/M_{i+1}$.
\end{Remark}
% The following example demonstrates that the above degree shift may not be zero.
% \begin{Example}
%     Take $\Lambda=\Lambda_0$, $\Gamma=\Aone[2]$, $\lambda=(4,1^5)$, $\lambda_1=(5,1^4)$. Let $\bi=\res(T^\lambda)$, then:
%     \[\deg T^{\lambda}=1\]
%     \[\deg T^{\lambda_1}=1\]
%     \[\deg \psi_1\psi_2\psi_3\psi_4e(\bi)=-1\]
%     Thus $-\deg T^{\lambda_i}+\deg T^{\lambda}+\deg \psi^{A_i}e(\bi)\neq 0$.
% \end{Example}

\begin{Remark}\label{rmk:skew-specht-filtration-hook}
    The Specht filtration in \autoref{thm-specht-filtration-hook} is generated by the Garnir relations in top‐to‐bottom order. By reversing this order, one obtains a skew Specht filtration of $M^\lambda$ whose factors are the skew Specht modules introduced in \cite{muth-graded-skew-specht}.
\end{Remark}
%%%%%%%%%%%%%%%%%%%%%%%%%%%%%%%%%%%%%%%%%%%%%%%%%%%%%%%%%%%%%%%%%%%%%%%%%%%%%%%%%%%%%%%%%%%%%%%%%%%%%%%%%%%%%%%%%%%%%%%%%%%%%%%%%%%%%%%%%%%%%%%%%%%%%%%%%%%%%%%%%%%%%%%%%%%%%%%%%%%%%%%%%%%%%%%%%%%%%%%%%%%%%%%%%%%%%%%%%%%%%%%%%%%%%%%%%%%%%%%%%%%%%%%%%%%%%%%%%%%%%%%%%%%%%%%%%%%%%%%%%%%%%%%%%%%%%%%%%%%%%%%%%%%%%%%%%%%%%%%%%%%%%%%%%%%%%%%%%%%%%%%%%%%%%%%%%%%%%%%%%%%%%%%%%%%%%%%%%%%%%%%%%%%%%%%%%%%%%%%%%%%%%%%%%%%%%%%%%%%%%%%%%%%%%%%%%%%%%%%%%%%%%%%%%%%%%%%%%%%%%%%%%%%%%%%%%%%%%%%%%%%%%%%%%%%%%%%%%%%%%%%%%%%%%%%%%%%%%%%%%%%%%%%%%%%%%%%%%%%%%%%%%%%%%%%%%%%%%%%%%%%%%%%%%%%%%%%%%%%%%%%%%%%%%%%%%%%%%%%%%%%%%%%%%%%%%%%%%%%%%%%%%%%%%%%%%%%%%%%%%%%%%%%%%%%%%%%%%%%%%%
\section{Two-Row Partition Case in Type \texorpdfstring{$A_\infty$}{A\_infinity}}\label{sec:specht-filtration-2row}
In this section, we fix a partition $\lambda=(k,r)$ with $k\geq r$ and prove that, in the linear quiver $A_\infty$ case (or in type $\Aone[e-1]$ with $e$ large enough), all Garnir relations of $\lambda$ are generated by the first one.

With the data above, the initial tableau $T^\lambda$ is the following:
\begin{center}
    \resizebox{\textwidth}{!}{
    \Tableau[xscale=2.5,dotted cols={3,4,7,8,11,12}]{1200{s}{s{+}1}00{r}{r{+}1}00{k{-}1}{k},{k{+}1}{k{+}2}00{k{+}s}{k{+}s{+}1}00{k{+}r}}
    }
\end{center}
Let $B_i:=(1,i)$ for $1\leq i\leq r$. Clearly they are all the Garnir nodes of $[\lambda]$. Set $\psi^{B_i}:=\psi^{G^{B_i}}$ for each $1\le i\le r$, and $\bi:=\bi^\lambda$.

\begin{Lemma}\label{lm:garnir-relation-2row}
    For $1\leq s\leq r$, we have:
    \begin{equation}
        \psi^{B_s}=(\psi_{2s-1}\cdots\psi_{k+s-2}\psi_{k+s-1})\cdots(\psi_{s+1}\cdots\psi_k\psi_{k+1})(\psi_s\cdots\psi_{k-1}\psi_{k})e(\bi)
    \end{equation}
\end{Lemma}
\begin{proof}
    If $s=1$, the permutation $\sigma_1\cdots \sigma_k$ transforms the initial tableau $T^{\lambda}$ into the following Garnir tableau $G^{B_1}$:
    \begin{center}
        \resizebox{\textwidth}{!}{
        \Tableau[xscale=2.5,dotted cols={4,5,6,7,8,9,12,13}]{234000000{r{+}1}{r{+}2}00{k}{k{+}1},{1}{k{+}2}{k{+}3}000000{k{+}r}}
        }
    \end{center}

   If $s>1$, the permutation $\sigma_s \cdots \sigma_{k-1} \sigma_k$ transforms the initial tableau into the following:
    \begin{center}
        \resizebox{\textwidth}{!}{
        \Tableau[xscale=2.5,dotted cols={4,5,8,9,12,13}]{12300{s{+}1}{s{+}2}00{r{+}1}{r{+}2}00{k}{k{+}1},{s}{k{+}2}{k{+}3}00{k{+}s}{k{+}s{+}1}00{k{+}r}}
        }
    \end{center}
    %(\psi_{s+t}\cdots\psi_{k+t-1}\psi_{k+t})\cdots
    Similarly, after applying $(\sigma_{s+1} \cdots \sigma_k \sigma_{k+1})(\sigma_s \sigma_{k-1} \cdots \sigma_k) e(\bi)$ to $T^\lambda$, we obtain the following tableau:
    \begin{center}
        \resizebox{\textwidth}{!}{
        \Tableau[xscale=2.5,dotted cols={4,5,8,9,12,13}]{12300{s{+}2}{s{+}3}00{r{+}2}{r{+}3}00{k{+}1}{k{+}2},{s}{s{+}1}{k{+}3}00{k{+}s}{k{+}s{+}1}00{k{+}r}}
        }
    \end{center}
    Continuing this procedure $s$ times, we obtain the Garnir tableau $G^{B_s}$:
    \begin{center}
        \resizebox{\textwidth}{!}{
        \Tableau[xscale=3,dotted cols={4,5,8,9,12,13}]{12300{2s}{2s{+}1}00{s{+}r}{s{+}r{+}1}00{k{+}s{-}1}{k{+}s},{s}{s{+}1}{s{+}2}00{2s{-}1}{k{+}s{+}1}00{k{+}r}}
        }
    \end{center}
    Let  
    \[
        w:=(\sigma_{2s-1}\cdots\sigma_{k+s-2}\sigma_{k+s-1})\cdots(\sigma_{s+1}\cdots\sigma_k\sigma_{k+1})(\sigma_s\cdots\sigma_{k-1}\sigma_{k})\in \Sym_{k+r}.
    \]
    We have shown that $wT^{\lambda}=G^{B_s}$. It is classical (see, for instance, \cite[Proposition 3.3]{mathas-iwahori-hecke}) that the map
    \[
    \{\text{row-standard } \lambda\text{-tableaux}\}\ \longrightarrow\ W^{\lambda}\subset \Sym_{k+r},\qquad T\longmapsto w^T,
    \]
    is a bijection onto the set $W^{\lambda}$ of minimal-length representatives of $\Sym_{k+r}/(\Sym_k\times \Sym_r)$, characterized by $w^T T^\lambda=T$. Therefore, since $wT^\lambda=G^{B_s}$, we must have $w=w^{G^{B_s}}$, and hence $w\in W^{\lambda}$. In particular, $w$ has minimal length in its coset and therefore the expression is reduced.
    
    Moreover, by \cite[Lemma 3.17]{kmr-universal-specht-type-A}, $w$ is \emph{fully commutative}, so any two reduced expressions for $w$ are related by commuting braid moves. Consequently, $\psi^w$ is independent of the choice of reduced expression, and the stated equality follows.
\end{proof}
\begin{Lemma}\label{lm:psi-s+1-kill-2s}
    For $1\leq s\leq r-1$, take the Garnir node $B_{s+1}$ and let $v$ be the standard cyclic generator of $M^{\lambda}$. Then
    \[\psi_{2s}\psi^{B_{s+1}}e(\bi)v=0\]
\end{Lemma}
\begin{proof}    
    If $s>1$, set $\psi':=(\psi_{2s-1}\cdots \psi_{k+s-2})\cdots(\psi_{s+1}\cdots\psi_k)$. By \autoref{lm:garnir-relation-2row}, we have:
    \[\psi^{B_{s+1}}=(\psi_{2s+1}\cdots\psi_{k+s})(\psi_{2s}\cdots\psi_{k+s-1})\psi'.\]
    Let $T_1:=(\sigma_{2s-1}\cdots \sigma_{k+s-2})\cdots(\sigma_{s+1}\cdots\sigma_k)T^{\lambda}$, then it is of the following form:
    \begin{center}
        \resizebox{\textwidth}{!}{
        \Tableau[xscale=3,dotted cols={3,4,9,10,12,13}]{1200{s{-}1}{s}{2s}{2s{+}1}00{s{+}r{-}1}00{k{+}s{-}1},{s{+}1}{s{+}2}00{2s{-}1}{k{+}s}{k{+}s{+}1}{k{+}s{+}2}00{k{+}r}}
        }
    \end{center}
    Let $T_2:=(\sigma_{t}\sigma_{t+1}\cdots\sigma_{k+s})(\sigma_{t-1}\cdots\sigma_{k+s-1})T_1$ for $2s+2\leq t< k+s$, then it is of the following form:
    \begin{center}
        \resizebox{\textwidth}{!}{\Tableau[xscale=3,dotted cols={3,8,11,13}]{120{s{-}1}{s}{2s}{2s{+}1}0{t{-}2}{t{+}1}0{s{+}r{-}1}0{k{+}s{+}1},{s{+}1}{s{+}2}0{2s{-}1}{t{-}1}{t}{k{+}s{+}2}0{2t{-}2{-}2s}{2t{-}1{-}2s}0{k{+}r}}
        }
    \end{center}
    It is easy to see $\res_{T_2}(t-2)\neq \res_{T_2}(t)$ and $\res_{T_1}(2s)\neq \res_{T_1}(2s+2)$ since the quiver is $A_\infty$ (or $\Aone[e-1]$ with $e\gg0$). We can apply \autoref{eq:braid_relation}:
    \begin{align*}
      \psi_{2s}\psi^{B_{s+1}}e(\bi)v&=\psi_{2s}(\psi_{2s+1}\psi_{2s+2}\cdots\psi_{k+s})(\psi_{2s}\psi_{2s+1}\cdots\psi_{k+s-1})\psi'e(\bi)v\\
        &=\bigl(\psi_{2s}\psi_{2s+1}\psi_{2s}\bigr)(\psi_{2s+2}\cdots\psi_{k+s})(\psi_{2s+1}\cdots\psi_{k+s-1})\psi'e(\bi)v\\
        &=\bigl(\psi_{2s+1}\psi_{2s}\psi_{2s+1}\bigr)(\psi_{2s+2}\cdots\psi_{k+s})(\psi_{2s+1}\cdots\psi_{k+s-1})\psi'e(\bi)v\\
        &=(\psi_{2s+1}\psi_{2s})\bigl(\psi_{2s+1}\psi_{2s+2}\psi_{2s+1}\bigr)(\psi_{2s+3}\cdots\psi_{k+s})(\psi_{2s+3}\cdots\psi_{k+s-1})\psi'e(\bi)v\\
        &=\cdots\\
        &=(\psi_{2s+1}\psi_{2s})(\psi_{2s+2}\psi_{2s+1})\cdots(\psi_{k+s-1}\psi_{k+s-2})\psi_{k+s-1}\psi_{k+s}\psi_{k+s-1}\psi'e(\bi)v\\
        &=(\psi_{2s+1}\psi_{2s})(\psi_{2s+2}\psi_{2s+1})\cdots(\psi_{k+s-1}\psi_{k+s-2})\psi_{k+s}\psi_{k+s-1}\psi_{k+s}\psi'e(\bi)v\\
        &=(\psi_{2s+1}\psi_{2s})(\psi_{2s+2}\psi_{2s+1})\cdots(\psi_{k+s-1}\psi_{k+s-2})\psi_{k+s}\psi_{k+s-1}\psi'\psi_{k+s}e(\bi)v=0.
    \end{align*}
    The second last equality holds because $\psi_{k+s}$ commutes with $\psi'$ and it kills $v$. 

    If $s=1$, by \autoref{lm:garnir-relation-2row}, $\psi^{B_2}=(\psi_3\cdots\psi_{k+1})(\psi_2\cdots\psi_{k})e(\bi)$. For $4\le i\le k+1$, the tableau $T:=(\sigma_{i}\cdots\sigma_{k+1})(\sigma_{i-1}\cdots\sigma_{k})T^{\lambda}$ is as follows:
    \begin{center}
        \resizebox{\textwidth}{!}{
        \Tableau[xscale=2.5,dotted cols={4,5,9,10,13,14}]{12300{i{-}2}{i{+}1}{i{+}2}00{r{+}2}{r{+}3}00{k{+}1}{k{+}2},{i{-}1}{i}{k{+}3}00{k{+}i{-}2}{k{+}i{-}1}{k{+}i}00{k{+}r}}
        }
    \end{center}
    In particular, $\res_{T}(i-2)\neq \res_{T}(i)$. We can apply \autoref{eq:braid_relation} and compute:
    \begin{align*}
        \psi_2\psi^{B_2}e(\bi)v&=\psi_2(\psi_3\cdots\psi_{k+1})(\psi_2\cdots\psi_{k})e(\bi)v\\
        &=(\psi_2\psi_3\psi_2)(\psi_4\cdots\psi_{k+1})(\psi_3\cdots\psi_{k})e(\bi)v\\
        &=(\psi_3\psi_2\psi_3)(\psi_4\cdots\psi_{k+1})(\psi_3\cdots\psi_{k})e(\bi)v\\
        &=(\psi_3\psi_2)(\psi_3\psi_4\psi_3)(\psi_5\cdots\psi_{k+1})(\psi_4\cdots\psi_{k})e(\bi)v\\
        &=(\psi_3\psi_2)(\psi_4\psi_3\psi_4)(\psi_5\cdots\psi_{k+1})(\psi_4\cdots\psi_{k})e(\bi)v\\
        &=\cdots\\
        &=(\psi_3\psi_2)(\psi_4\psi_3)\cdots (\psi_{k}\psi_{k-1})\psi_k\psi_{k+1}\psi_ke(\bi)v\\
        &=(\psi_3\psi_2)(\psi_4\psi_3)\cdots (\psi_{k}\psi_{k-1})\psi_{k+1}\psi_k \psi_{k+1}e(\bi)v\\
        &=0.
    \end{align*}

    % Observe that $G^{B_{s+1}}$ is the following tableau:
    % \begin{center}
    %     \Tableau[xscale=2.1,dotted cols={4,5,10,11,14,15}]{12300{s{-}1}{s}{2s{+}2}{2s{+}3}00{s{+}r{+}1}{s{+}r{+}2}00{k{+}s}{k{+}s{+}1},{s{+}1}{s{+}2}{s{+}3}00{2s{-}1}{2s}{2s{+}1}{k{+}s{+}2}00{k{+}r}}
    % \end{center}
    % and the tableau $T':=(\sigma_s\sigma_{s+1}\cdots\sigma_{2s-2}\sigma_{2s-1}\sigma_{2s-2}\cdots\sigma_{s+1}\sigma_s)G^{B_{s+1}}$ is the following:
    % \begin{center}
    %     \Tableau[xscale=2.1,dotted cols={4,5,10,11,14,15}]{12300{s{-}1}{2s{-}1}{2s}{2s{+}3}00{s{+}r{+}1}{s{+}r{+}2}00{k{+}s}{k{+}s{+}1},{s{+}1}{s{+}2}{s{+}3}00{2s{-}1}{s}{2s{+}1}{k{+}s{+}2}00{k{+}r}}
    % \end{center}
    % As $\res_{T'}()$
    % Let $\bj:=\res(G^{B_{s+1}})$, we know $\psi_{2s}\psi^{B_{s+1}}e(\bi)=\psi_{2s}e(\bj)\psi^{B_{s+1}}=e(\sigma_{2s}\bj)\psi_{2s}\psi^{B_{s+1}}$. 
\end{proof}
We are now in a position to prove the main result of this section:
\begin{Theorem}\label{thm:2row-garnir-connection}
    Let $v$ be the standard cyclic generator of $M^\lambda$, for $1\leq s\leq r-1$, we have:
    \begin{equation}
        (\psi_{s}\psi_{s+1}\cdots\psi_{k+s})\psi^{B_s}e(\bi)v=-\psi^{B_{s+1}}e(\bi)v
    \end{equation}
\end{Theorem}
\begin{proof}
    We compute as follows:
    
    $(\psi_{s}\psi_{s+1}\cdots\psi_{k+s})\psi^{B_s}e(\bi)v$
    \begin{align*}
        &=(\psi_{s}\psi_{s+1}\cdots\psi_{k+s})(\psi_{2s-1}\cdots\psi_{k+s-1})\cdots(\psi_{s+1}\cdots\psi_{k+1})(\psi_s\cdots\psi_{k})e(\bi)v\\
        %%%%%%%%
        &=(\psi_{s}\psi_{s+1}\cdots\psi_{k+s})(\psi_{2s-1}\psi_{2s-2}\cdots\psi_s)(\psi_{2s}\cdots\psi_{k+s-1})\cdots(\psi_{s+2}\cdots\psi_{k+1})(\psi_{s+1}\cdots\psi_{k})e(\bi)v\\
        %%%%%%%%
        &=(\psi_{s}\psi_{s+1}\cdots\psi_{2s-1}\psi_{2s}\psi_{2s+1}\cdots\psi_{k+s})(\psi_{2s-1}\psi_{2s-2}\cdots\psi_s)(\psi_{2s}\cdots\psi_{k+s-1})\cdots(\psi_{s+2}\cdots\psi_{k+1})(\psi_{s+1}\cdots\psi_{k})e(\bi)v\\
        %%%%%%%%
        &=(\psi_{s}\psi_{s+1}\cdots\psi_{2s-1}\psi_{2s}\psi_{2s-1}\psi_{2s-2}\cdots\psi_s)(\psi_{2s+1}\cdots\psi_{k+s}))(\psi_{2s}\cdots\psi_{k+s-1})\cdots(\psi_{s+2}\cdots\psi_{k+1})(\psi_{s+1}\cdots\psi_{k})e(\bi)v\\
        %%%%%%%%
        &=(\psi_{s}\psi_{s+1}\cdots\psi_{2s-1}\psi_{2s}\psi_{2s-1}\psi_{2s-2}\cdots\psi_s)\psi^{B_{s+1}}e(\bi)v.\\
    \end{align*}
    
    Notice $T^{2s-2}:=(\sigma_{2s-2}\cdots\sigma_s)G^{B_{s+1}}$ is the following tableau:
    \begin{center}
        \resizebox{\textwidth}{!}{\Tableau[xscale=2.5,dotted cols={4,9,12}]{1230{s{-}1}{2s{-}1}{2s{+}2}{2s{+}3}0{s{+}r{+}1}{s{+}r{+}2}0{k{+}s}{k{+}s{+}1},{s}{s{+}1}{s{+}2}0{2s{-}2}{2s}{2s{+}1}{k{+}s{+}2}0{k{+}r}}
        }
    \end{center}
    Since $\res_{T^{2s-2}}(2s-1)=\res_{T^{2s-2}}(2s+1)\leftarrow \res_{T^{2s-2}}(2s)$, we get:
    \begin{align*}
        &\hspace{5mm}(\psi_{s}\psi_{s+1}\cdots\psi_{2s-1}\psi_{2s}\psi_{2s-1}\psi_{2s-2}\cdots\psi_s)\psi^{B_{s+1}}e(\bi)v\\
    %%%%%%%
        &=(\psi_{s}\psi_{s+1}\cdots\psi_{2s-2}\bigl(\psi_{2s}\psi_{2s-1}\psi_{2s}-1\bigr)\psi_{2s-2}\cdots\psi_s)\psi^{B_{s+1}}e(\bi)v\\
    %%%%%%%
        &=(\psi_{s}\psi_{s+1}\cdots\psi_{2s-2}\bigl(\psi_{2s}\psi_{2s-1}\psi_{2s}\bigr)\psi_{2s-2}\cdots\psi_s)\psi^{B_{s+1}}e(\bi)v-(\psi_{s}\psi_{s+1}\cdots\psi_{2s-2}\psi_{2s-2}\cdots\psi_s)\psi^{B_{s+1}}e(\bi)v.
    \end{align*}
    It is clear the first term vanishes by \autoref{lm:psi-s+1-kill-2s} as $\psi_{2s}$ commutes with the remaining $\psi_{2s-2}\cdots\psi_s$. 

    Let $T^r:=(\sigma_{r}\cdots\sigma_s)G^{B_{s+1}}$ for $s< r< 2s-2$, which is of the following form:
    \begin{center}
        \resizebox{\textwidth}{!}{\Tableau[xscale=2.5,dotted cols={3,6,9,12}]{120{r{-}s{+}1}{r{-}s{+}2}0{r{+}1}{2s{+}2}0{s{+}r{+}1}{s{+}r{+}2}0{k{+s{+}1}},s{s{+}1}0r{r{+2}}0{2s}{2s{+}1}0{k{+}r}}
        }
    \end{center}
    As $\res_{T^r}(r+1)$ and $\res_{T^r}(r+2)$ are not adjacent or equal, apply \autoref{eq:psi_square}:
    \begin{align*}
        (\psi_{s}\psi_{s+1}\cdots\psi_{2s-2}\psi_{2s-2}\cdots\psi_s)\psi^{B_{s+1}}e(\bi)v
        &=(\psi_{s}\psi_{s+1}\cdots\psi_{2s-3}\psi_{2s-3}\cdots\psi_s)\psi^{B_{s+1}}e(\bi)v\\
        &=\cdots\\
        &=\psi^{B_{s+1}}e(\bi)v.
    \end{align*}
\end{proof}
\begin{Corollary}\label{cor:garnir-generator-2row}
    The submodule of $M^\lambda$ generated by
    \[
        \{ \psi^{B_1} e(\bi) v, \psi^{B_2} e(\bi) v, \ldots, \psi^{B_r} e(\bi) v \}
    \]
    is cyclic, with generator $\psi^{B_1} e(\bi) v$.
\end{Corollary}
%%%%%%
\begin{Corollary}\label{cor:garnir-relations-simplify}
    Fix quiver $A_\infty$ or $\Aone[e-1]$ with $e\gg0$. Fix $\Lambda\in P^+$ a fundamental weight and $\alpha\in Q^+$ a positive root. Let $\lambda=(\lambda_1,\cdots,\lambda_r)\in \Par[\Lambda]_{\alpha}$ and $\bi:=\bi^\lambda$, and form the permutation module $M^\lambda$ with standard cyclic generator $v$, then the submodule $M_1$ of $M^\lambda$ generated by the Garnir relations has the following form:
    \[
        M_1=R_{\alpha}\{ \psi^{A_1}e(\bi)v,\cdots,\psi^{A_{r-1}}e(\bi)v\},
    \]
    where $\{A_i=(i,1)|1\leq i\leq r-1\}$ is the set of Garnir nodes in the first column of $[\lambda]$. 
\end{Corollary}
\begin{proof}
    Each Garnir relation only involves two rows, hence we only need to show any Garnir relation relating row $i$ and $i+1$ is generated by the Garnir relation corresponding to the Garnir node $A_i$, which is just \autoref{cor:garnir-generator-2row}.
\end{proof}
%%%%%%
We are now ready to construct a Specht filtration of $M^\lambda$. As a first step, we compare the dimensions, analogous to \autoref{lm:equality-dimension-specht-permutation}.
\begin{Lemma}\label{lm:dimension-equality-2row}
    Let $\lambda_0=\lambda$ and $\lambda_1=(k+1|r-1)$. Then $\dim M^\lambda=\dim S^{\lambda_0}+\dim S^{\lambda_1}$.
\end{Lemma}
\begin{proof}
    The dimension of $M^\lambda$ is \[\frac{(k+r)!}{k!r!},\] and the dimension of $S^{\lambda_0}$ is \[\frac{(k+r)!}{\bigl((k+1)\cdots (k-r+2)\bigr)(k-r)!r!}=\frac{(k+r)!(k-r+1)}{(k+1)!r!}.\]
    The dimension of $S^{\lambda_1}$ is:
    \[\binom{k+r}{k+1}=\frac{(k+r)!}{(k+1)!(r-1)!}=\frac{(k+r)!r}{(k+1)!r!}.\]
    It is immediate to get the desired equality.
\end{proof}

\begin{Theorem}\label{thm:specht-filtration-2row}
    Suppose the quiver is $A_\infty$ or $\Aone[e-1]$ with $e \gg 0$. Fix $\alpha\in Q^+$ and $\Lambda_i\in P^+$. Let $\lambda = (k, r)$ with $k \geq r$ be a two-row partition such that $\alpha_{\lambda}=\alpha$. Let $v$ be the standard cyclic generator of the permutation module $M^\lambda$, and let $\bi := \bi^\lambda$. Then $M^\lambda$ admits a Specht filtration:
    \[
        M^\lambda = M_0 \supsetneq M_1 \supsetneq 0,
    \]
    where $M_1 = R_\alpha \cdot \psi^{B_1} e(\bi) v$, and we have:
    \[
        M_0 / M_1 \cong S^{\lambda}, \qquad M_1 \cong S^{\lambda_1},
    \]
    where $\lambda_1 = (k+1 \mid r-1)$ and $S^{\lambda_1}$ is the Specht module over $R_{\alpha}^{\Lambda}$ with $\Lambda$ determined by the charge $\kappa=({i-1},i)$.
\end{Theorem}

\begin{proof}
    By \autoref{cor:garnir-generator-2row}, the only part that remains to be proven is the isomorphism between $M_1$ and $S^{\lambda_1}$. Using essentially the same argument as in \autoref{lm-proof-i=1}, we can show that the cyclic generator $\psi^{B_1} e(\bi) v$ satisfies all the defining relations of the standard cyclic generator $w$ of $S^{\lambda_1}$. Hence, there exists a surjective homomorphism from $S^{\lambda_1}$ onto $M_1$. The conclusion then follows from \autoref{lm:dimension-equality-2row}.
\end{proof}

At the end of this section, we record the following result:

\begin{Lemma}\label{lm:garnir-connection-2}
    Let $v$ be the standard cyclic generator of $M^\lambda$. For $2 \leq s \leq r$, we have:
    \[
        (\psi_{k+s-1} \cdots \psi_{s} \psi_{s-1})\, \psi^{B_s} e(\bi) v = -\psi^{B_{s-1}} e(\bi) v.
    \]
\end{Lemma}
%%%%%%%%%%%%%%%%%
\begin{proof}
    The proof is analogous to that of \autoref{thm:2row-garnir-connection}. We briefly state the key procedures:
    \begin{align*}
        &(\psi_{k+s-1} \cdots \psi_{s} \psi_{s-1})\, \psi^{B_s} e(\bi) v\\
        &= (\psi_{k+s-1} \cdots \psi_{s} \psi_{s-1})\, (\psi_{2s-1}\cdots\psi_{k+s-2}\psi_{k+s-1})\cdots(\psi_{s+1}\cdots\psi_k\psi_{k+1})(\psi_s\cdots\psi_{k-1}\psi_{k})e(\bi) v\\
        &=(\psi_{k+s-1}\cdots\psi_{2s-2})(\psi_{2s-1}\cdots\psi_{k+s-2}\psi_{k+s-1})\psi^{B_{s-1}}e(\bi)v\\
        &=(\psi_{k+s-1}\cdots\psi_{2s-2}\psi_{2s-1}\psi_{2s-2}\cdots\psi_{k+s-2}\psi_{k+s-1})\psi^{B_{s-1}}e(\bi)v\\
        &=\bigl(\psi_{k+s-1}\cdots\psi_{2s-3}(\psi_{2s-1}\psi_{2s-2}\psi_{2s-1}-1)\psi_{2s-3}\cdots\psi_{k+s-2}\psi_{k+s-1}\bigr)\psi^{B_{s-1}}e(\bi)v.
    \end{align*}
    Then we need to prove $\psi_{2s-2}\psi^{B_{s-1}}e(\bi)v=0$ and the conclusion follows.
\end{proof}

\autoref{lm:garnir-connection-2}, together with \autoref{thm:2row-garnir-connection}, shows that the choice of Garnir relation in \autoref{cor:garnir-generator-2row} between any two adjacent rows is be arbitrary. In other words, for each pair of adjacent rows $i$ and $i+1$, we may choose any Garnir node $(i,j) \in [\lambda]$ and use the corresponding Garnir relation. This single relation suffices to generate all the Garnir relations between the two rows in $[\lambda]$. 
%%%%%%%%%%%%%%%%%%%%%%%%%%%%%%%%%%%%%%%%%%%%%%%%%%%%%%%%%%%%%%%%%%%%%%%%%%%%%%%%%%%%%%%%%%%%%%%%%%%%%%%%%%%%%%%%%%%%%%%%%%%%%%%%%%%%%%%%%%%%%%%%%%%%%%%%%%%%%%%%%%%%%%%%%%%%%%%%%%%%%%%%%%%%%%%%%%%%%%%%%%%%%%%%%%%%%%%%%%%%%%%%%%%%%%%%%%%%%%%%%%%%%%%%%%%%%%%%%%%%%%%%%%%%%%%%%%%%%%%%%%%%%%%%%%%%%%%%%%%%%%%%%%%%%%%%%%%%%%%%%%%%%%%%%%%%%%%%%%%%%%%%%%%%%%%%%%%%%%%%%%%%%%%%%%%%%%%%%%%%%%%%%%%%%%%%%%%%%%%%%%%%%%%%%%%%%%%%%%%%%%%%%%%%%%%%%%%%%%%%%%%%%%%%%%%%%%%%%%%%%%%%%%%%%%%%%%%%%%%%%%%%%%%%%%%%%%%%%%%%%%%%%%%%%%%%%%%%%%%%%%%%%%%%%%%%%%%%%%%%%%%%%%%%%%%%%%%%%%%%%%%%%%%%%%%%%%%%%%%%%%%%%%%%%%%%%%%%%%%%%%%%%%%%%%%%%%%%%%%%%%%%%%%%%%%%%%%%%%%%%%%%%%%%%%%%%%%%%%%%%%
\section{General Partition Case in Type \texorpdfstring{$A_\infty$}{A\_infinity}}\label{sec:specht-filtration-general-partition}
One might expect that \autoref{thm:specht-filtration-2row} naturally extends to arbitrary partitions, particularly given \autoref{cor:garnir-relations-simplify}. However, this extension encounters a fundamental obstacle: the dimension equality established in \autoref{lm:dimension-equality-2row} fails to hold for general partitions. In fact, the surjective homomorphism from $S^{\lambda_i}$ to $M_i/M_{i+1}$ is not an isomorphism in general. Instead, we construct a finite Specht resolution of $M_i/M_{i+1}$.
\begin{Example}
    Fix $\Lambda_0$ and consider the partition $\lambda=\lambda_0=(5,5,4,2,2)$.
    \begin{center}
        \Tableau{01234,{{-}1}0123,{{-}2}{{-}1}01,{{-}3}{{-}2},{{-}4}{{-}3}}
    \end{center}
    We have $\dim M^\lambda=4631346720$ and $\dim S^\lambda=4594590$.
    
    The `filtration' generated by the Garnir relations contains the following Specht modules:

    $\lambda_1=\big(4|(6,4,2,2)\big)$:
    \begin{center}
        \Multitableau{0123|{{-}1}01234,{{-}2}{{-}1}01,{{-}3}{{-}2},{{-}4}{{-}3}}
    \end{center}
    $\lambda_2=\big(5|3|(6,2,2)\big)$:
    \begin{center}
        \Multitableau{01234|{{-}1}01|{{-}2}{{-}1}0123,{{-}3}{{-}2},{{-}4}{{-}3}}
    \end{center}
    $\lambda_3=\big(5|5|1|(5,2)\big)$:
    \begin{center}
        \Multitableau{01234 | {{-}1}0123 | {{-2}} | {{-}3}{{-}2}{{-}1}01,{{-}4}{{-}3}}
    \end{center}
    $\lambda_4=\big(5|5|4|1|3\big)$:
    \begin{center}
        \Multitableau{01234|{{-}1}0123|{{-}2}{{-}1}01|{{-3}}|{{-}4}{{-}3}{{-}2}}
    \end{center}
    By computation, we get:
    \[\dim S^{\lambda_1}=128648520,\]
    \[\dim S^{\lambda_2}=551350800,\]
    \[\dim S^{\lambda_3}=1235025792,\]
    \[\dim S^{\lambda_4}=3087564480.\]
    However, unlike the case of hook or two-row partitions, this time we have:
    \[
        \sum\limits_{0 \leq i \leq 4} \dim S^{\lambda_i} \neq \dim M^\lambda.
    \]
    The problem arises because, the surjection $S^{\lambda_i} \twoheadrightarrow M_i / M_{i+1}$ is not an isomorphism in general. Indeed, the kernels of these maps are themselves Specht modules corresponding to the following partitions:\\    
    $\mu_1=\big(2|(6,6,2,2)\big)$:
    \begin{center}
        \Multitableau{01|{{-}1}01234,{{-}2}{{-}1}0123,{{-}3}{{-}2},{{-}4}{{-}3}}\\
    \end{center}
    $\mu_2=\big(5|\emptyset|(6,5,2)\big)$:
    \begin{center}
        \Multitableau[empty=$\emptyset$]{01234||{{-}2}{{-}1}0123,{{-}3}{{-}2}{{-}1}01,{{-}4}{{-}3}}\\
    \end{center}
    $\mu_3=\big(5|5|\emptyset|(5,3)\big)$:
    \begin{center}
        \Multitableau[empty=$\emptyset$]{01234|{{-}1}0123||{{-}3}{{-}2}{{-}1}01,{{-}4}{{-}3}{{-}2}}\\
    \end{center}
    and the dimensions are:
    $\dim S^{\mu_1}=22972950$, $\dim S^{\mu_2}=44108064$ and  $\dim S^{\mu_3}=308756448$.
    
    It is not hard to see that:
    \[
        \sum\limits_{0 \leq i \leq 4} \dim S^{\lambda_i} - \sum\limits_{1 \leq i \leq 3} \dim S^{\mu_i} = \dim M^\lambda.
    \]
\end{Example}
Our main result in this section is the following:
\begin{Theorem}\label{thm:specht-filtration-general-partition}
    Suppose the quiver is $A_\infty$ or $\Aone[e-1]$ with $e \gg 0$. Fix $\alpha \in Q^+$ and let $\Lambda := \Lambda_x$ be a fundamental weight. Take $\lambda = (\lambda_1, \cdots, \lambda_r)\in \Par[\Lambda]_{\alpha}$. Set
    \[
        k_i=k_i(\lambda) :=
            \begin{cases}
                \max\{1 \leq j \leq r - i \mid \lambda_{i + j} - j \geq 0\}, & \text{if } 1 \leq i \leq r - 1, \\
                1, & \text{if } i = 0.
            \end{cases}
    \]
    Let $v$ be the standard cyclic generator of the permutation module $M^\lambda$, and let $\bi := \bi^\lambda$. Define $B_i := (i, \lambda_{i+1}) \in [\lambda]$ for $1 \leq i \leq r - 1$. Then $M^\lambda$ admits a generalized Specht filtration in the following sense:
    \[
        M^\lambda = M_0 \supsetneq M_1 \supsetneq \cdots \supsetneq M_{r-1} \supsetneq M_r = 0
    \]
    such that for each $0 \leq i \leq r - 1$, there exists an exact sequence of $R_\alpha$-modules:
    \begin{equation}\label{eq:specht-resolution-thm}
        0 \to S^{\mu_{i, k_i}} \xrightarrow{\upphi_{i, k_i}} S^{\mu_{i, k_i-1}}\cdots \xrightarrow{\upphi_{i,2}} S^{\mu_{i,1}} \xrightarrow{\upphi_{i,1}} M_i / M_{i+1} \to 0,
    \end{equation}
    where
    \[
        M_i = R_\alpha \{ \psi^{B_i} v, \cdots, \psi^{B_{r-1}} v \}, 1\leq i\leq r-1. 
    \]
    The multipartitions $\mu_{i,j}$ is given by
    \[
        \mu_{i,j} = \big(\lambda_1 | \cdots | \lambda_{i-1} | \lambda_{i+j} - j | (\lambda_i + 1, \cdots, \lambda_{i + j - 1} + 1, \lambda_{i + j + 1}, \cdots, \lambda_r)\big), \quad 1 \leq i \leq r - 1,
    \]
    and
    \[
        \mu_{0,1} = \lambda.
    \]
    Here, $S^{\mu_{i,j}}$ is the Specht module associated with $\mu_{i,j}$ over the cyclotomic KLR algebra $R^{\Lambda(i)}_\alpha$, where
    $\Lambda(i)$ is determined by the charge $\kappa=(x, {x - 1}, \cdots, {x - i+1},{x - i})$.
    %and the charge $\kappa(i)=(\kappa_1,\cdots,\kappa_)$ is given by
    % Moreover, for each $0 \leq i \leq r - 1$ and $1 \leq j \leq k_i$, let $u_{i,j}$ denote the standard cyclic generator of $S^{\mu_{i,j}}$. Define $n_i := \sum\limits_{1 \leq j \leq i} \lambda_j$ for $1 \leq i \leq r$ and $d_i:=\lambda_i-\lambda_{i+1}+1$ for admissible $i$. Then the $R_\alpha$-homomorphism is given by:
    % \[
    %     \upphi_{i,j}(u_{i,j}) = u_{i,j-1}, \quad 1 < j \leq k_i, \quad \text{and} \quad \upphi_{i,1}(u_{i,1}) = \psi^{B_i} v+M_{i+1}.
    % \]
\end{Theorem}
The resolution for each $M_i/M_{i+1}$ is called a \textbf{Specht resolution} and the filtration is called a \textbf{generalized Specht filtration}. 

%The proof begins by establishing a surjective homogeneous homomorphism $\upphi_{i,1}: S^{\mu_{i,1}} \to M_i/M_{i+1}$ that maps $u_{i,1}$ to $\psi^{B_i}v+M_{i+1}$. This procedure is analogous to that in \autoref{lm-proof-i=1} and \autoref{lm-proof-r>i>1} and involves checking that $\psi^{B_i}v+M_{i+1}$ satisfies all the defining relations for $u_{i,1}$. Similarly, we verify that each $\upphi_{i,j}$ is a well-defined $R_\alpha$-homomorphism for admissible $j$, and that $\upphi_{i,k_i}$ is injective. Furthermore, we must check that $\upphi_{i,j}\circ \upphi_{i,j-1}=0$ for each $j$ such that $2\leq j\leq k_i$, which establishes that it forms a complex.

The following proof relies on results established later in this section. We present it first because it offers greater clarity.
\begin{proof}[Proof of \autoref{thm:specht-filtration-general-partition}]
The proof proceeds by induction on the length $r$ of the partition $\lambda=(\lambda_1,\dots,\lambda_r)$. The base cases, $r=1$ is trivial and $r=2$ was established in \autoref{thm:specht-filtration-2row}. Assume $r\geq 3$, and let $\nu=(\lambda_2,\dots,\lambda_r)$. By the induction hypothesis, assume that \autoref{thm:specht-filtration-general-partition} holds for any partition of length less than or equal to $r-1$. In particular, $M^\nu$ possesses the desired generalized Specht filtration:
    \begin{equation}\label{eq:generalized-filtration}
        M^\nu = N_0 \supsetneq N_1 \supsetneq \dots \supsetneq N_{r-2} \supsetneq N_{r-1} = 0
    \end{equation}
Let $\beta=\alpha_{\lambda_1}$ and define $F:=\text{Ind}^{R_\alpha}_{R_\beta\otimes R_{\alpha-\beta}}$. According to \cite[Proposition 2.16]{khovanovlauda-klr-1}, $F$ is an exact functor. Let $S^\beta=L_\beta$ be the one-dimensional Specht module associated with $(\lambda_1)$ over $R^{\Lambda_x}_\beta$. The functor $L_\beta\boxtimes- := L_\beta\otimes_\bk-$ is also an exact functor (since $L_\beta$ is free over $\bk$), mapping the category of finite-dimensional $R_{\alpha-\beta}$-modules to the category of finite-dimensional $R_{\beta}\otimes R_{\alpha-\beta}$-modules.
For modules $D_1$ over $R_\beta$ and $D_2$ over $R_{\alpha-\beta}$, respectively, define $D_1\circ D_2 := F(D_1\boxtimes D_2)$. In particular, if $D_1=S^\beta$, we consider the module $S^{\beta}\circ D_2$. By \autoref{thm:kmr-high-level-specht}, if $\nu'$ is a partition such that $\alpha_{\nu'}=\alpha-\beta$, then $S^\beta\circ S^{\nu'}\cong S^{\mu'}$, where $\mu'=(\lambda_1|\nu')$.

Hence, the generalized Specht filtration \autoref{eq:generalized-filtration} of $M^\nu$ yields the following sequence of $R_\beta\otimes R_{\alpha-\beta}$-modules:
    \[
        S^\beta\boxtimes M^\nu = S^\beta\boxtimes N_0 \supsetneq S^\beta\boxtimes N_1 \supsetneq \cdots \supsetneq S^\beta\boxtimes N_{r-2} \supsetneq S^\beta\boxtimes N_{r-1} = 0.
    \]
Applying the exact functor $F$, we obtain the following filtration of $R_\alpha$-modules:
    \begin{equation}\label{eq:filtration-r-1-circ}
        S^\beta\circ M^\nu = S^\beta\circ N_0 \supsetneq S^\beta\circ N_1 \supsetneq \cdots \supsetneq S^\beta\circ N_{r-2} \supsetneq S^\beta\circ N_{r-1} = 0.
    \end{equation}
By \autoref{lm:from-Ni-to-Mi+1}, $M_{i+1}\cong S^\beta\circ N_i$ for $1\leq i\leq r-2$, and $M^\lambda = S^\beta\circ M^\nu$.
This filtration \autoref{eq:filtration-r-1-circ} can therefore be rewritten as:
    \[
        M^\lambda = M_0 \supsetneq M_1 \supsetneq \cdots \supsetneq M_{r-1} \supsetneq M_{r} = 0.
    \]
Suppose that for each $i$, the Specht resolution of the quotient module $N_i/N_{i+1}$ is given by:
    \begin{equation}\label{eq:specht-resolution-Ni}
        0\to S^{\nu_{i,t_i}}\to \cdots \to S^{\nu_{i,1}}\to N_i/N_{i+1}\to 0
    \end{equation}
where $t_i=k_i(\nu)$. Let $k_i=k_{i}(\lambda)$. By our construction in \autoref{thm:specht-filtration-general-partition}, it is clear that $t_i=k_{i+1}$ and $\mu_{i+1,j}=(\lambda_1|\nu_{i,j})$.
Applying the exact functor $F(S^{\beta}\boxtimes-) = S^\beta\circ-$ to this resolution \autoref{eq:specht-resolution-Ni} yields:
\begin{center}
    \begin{tikzcd}
        0 \arrow[r] &  S^{\beta}\circ S^{\nu_{i,t_i}} \arrow[r] \arrow[d,"\cong"]& \cdots \arrow[r] & S^{\beta}\circ S^{\nu_{i,1}} \arrow[r] \arrow[d,"\cong"] & S^{\beta}\circ (N_i/N_{i+1}) \arrow[r] \arrow[d,"\cong"]& 0 \\
        %%%%%%
        0 \arrow[r] & S^{\mu_{i+1,k_{i+1}}} \arrow[r] & \cdots \arrow[r] & S^{\mu_{i+1,1}} \arrow[r] & M_{i+1}/M_{i+2} \arrow[r] & 0
    \end{tikzcd}
\end{center}
For $2\leq i\leq r-1$, this gives to the desired Specht resolution for $M_i/M_{i+1}$ from \autoref{thm:specht-filtration-general-partition}. The remaining task is to demonstrate that for the submodule $M_1$ (where $M_1 \subsetneq M_0$, and $M_1$ is generated by all Garnir relations), the quotient $M_1/M_2$ admits the following desired Specht resolution:
    \begin{equation}\label{eq:specht-resolution-M1}
        0\to S^{\mu_{1,k_1}}\xrightarrow{\upphi_{1,k_1}}\cdots\to S^{\mu_{1,1}}\xrightarrow{\upphi_{1,1}} M_1/M_2\to 0
    \end{equation}
We construct this resolution by applying results from \cite{humathas-quiver-schur-linear-quiver}. First, all the $2$-partitions $\mu_{1,j}$ ($1 \leq j \leq k_1$) are Kleshchev by \autoref{lm:mu_j-is-kleshchev-partitions}. Thus, each $S^{\mu_{1,j}}$ has a unique irreducible head $D^{\mu_{1,j}}$. 

Let $d_{\lambda,\mu}$ be the decomposition number $[S^\lambda : D^\mu]$, where $\mu$ is a Kleshchev multipartition. By \autoref{cor:decomposition-number-dominance}, we have
\[
d_{\mu_{1,j},\mu} = 
\begin{cases}
  1 & \text{if } \mu = \mu_{1,j} \text{ or } \mu_{1,j+1}, \\
  0 & \text{otherwise}.
\end{cases}
\]

This implies that each $S^{\mu_{1,j}}$ ($1 \leq j < k_1$) has a composition series:
\[
0 \subsetneq D^{\mu_{1,j+1}} \subsetneq S^{\mu_{1,j}},
\]
such that
\[
S^{\mu_{1,j}} / D^{\mu_{1,j+1}} \cong D^{\mu_{1,j}},
\]
with $S^{\mu_{1,k_1}} \cong D^{\mu_{1,k_1}}$ itself irreducible.

Using this structure, the resolution can now be constructed explicitly. The map $\upphi_{1,k_1}$ is the canonical embedding of $S^{\mu_{1,k_1}}$ into the submodule $D^{\mu_{1,k_1}}$ of $S^{\mu_{1,k_1-1}}$. For $2 \leq j \leq k_1 - 1$, define the maps $\upphi_{1,j}$ by sending the submodule $D^{\mu_{1,j+1}}$ to zero, thereby inducing:
\[
S^{\mu_{1,j}} / D^{\mu_{1,j+1}} \cong D^{\mu_{1,j}} \hookrightarrow S^{\mu_{1,j-1}}.
\]
The map $\upphi_{1,1}$ is constructed similarly to \autoref{thm:specht-filtration-2row}, by sending the standard cyclic generator $u_{1,1}$ of $S^{\mu_{1,1}}$ to the element $\psi^{B_1} v + M_2$. It is routine (and analogous to the arguments in \autoref{lm-proof-i=1} and \autoref{lm-proof-i=r}) to verify that $\psi^{B_1} v + M_2$ satisfies all defining relations for $u_{1,1}$. Hence, the map $\upphi_{1,1}$ is a surjective homogeneous homomorphism. The composition series ensures that $S^{\mu_{1,1}}$ has exactly one proper non-trivial submodule $D^{\mu_{1,2}}$. By \autoref{lm:kernel-not-zero} the kernel of $\upphi_{1,1}$ is nonzero if $k_1 \neq 1$.  Therefore, they coincide. If $k_1=1$, then $S^{\mu_{1,1}}$ is simple and isomorphic to $M_1/M_2$. 

Thus, we have constructed an exact sequence given by \autoref{eq:specht-resolution-M1}, completing the proof by induction.
\end{proof}

\begin{Corollary}\label{cor:description-factor}
    Assume the same conditions as in \autoref{thm:specht-filtration-general-partition}, and define the modules $M_i$ as in the generalized Specht filtration there.
    For each $1\leq i\leq r-1$, set $\nu_i:=\big(\lambda_{i+1} - 1 | (\lambda_i + 1, \lambda_{i + 2}, \cdots, \lambda_r)\big)$.
    Then $M_i/M_{i+1}\cong S^{(\lambda_1)}\circ \cdots\circ S^{(\lambda_{i-1})}\circ D^{\nu_i}$. In particular, $M_{r-1}\cong S^{\mu_{r-1,1}}$ and $M_1/M_2\cong D^{\mu_{1,1}}$. 
\end{Corollary}
\begin{proof}
    By the proof of \autoref{thm:specht-filtration-general-partition}, the statement is true for $i=1$ since $M_1/M_2\cong D^{\mu_{1,1}}$ and $\nu_1=\mu_{1,1}$. Then the general case follows by induction.
\end{proof}
%%%%%%%%%%
\autoref{cor:description-factor} illustrates how far our generalized Specht filtration deviates from an actual Specht filtration: rather than each factor $M_i/M_{i+1}$ being isomorphic to a Specht module, it is isomorphic to an ``almost-Specht" module. We use the term ``almost-Specht" to emphasize that $S^{(\lambda_1)} \circ \cdots \circ S^{(\lambda_{i-1})} \circ S^{\nu_i}$ is a Specht module, and the difference lies only in the final term.
\begin{Corollary}\label{cor:dimension-formula}
    Assume the same conditions as in \autoref{thm:specht-filtration-general-partition}, and define the partitions $\mu_{1,j}$ as in the generalized Specht filtration there, for each $1\leq j\leq k_1$. Then:
    \[\dim D^{\mu_{1,j}}=\sum\limits_{j\leq s\leq k_1}(-1)^{s-j}\dim S^{\mu_{1,s}}\]
\end{Corollary}
\begin{proof}
    In the proof of \autoref{thm:specht-filtration-general-partition}, we observe that
    \[
        \dim S^{\mu_{1,j}} = \dim D^{\mu_{1,j}} + \dim D^{\mu_{1,j+1}} \quad \text{for } j \neq k_1(\lambda),
    \]
    and
    \[
        \dim S^{\mu_{1,k_1(\lambda)}} = \dim D^{\mu_{1,k_1(\lambda)}}.
    \]
    The stated equality then follows immediately.
\end{proof}
%%%%%%%%%%
We remind the reader that the dimensions of simple modules in the representation theory of KLR algebras are generally difficult to compute. In our setting, however, the dimensions of the Specht modules $S^{\mu_{1,j}}$ can be readily determined via the hook length formula. Thus, the last corollary provides a convenient method for computing the dimensions of the simple modules $D^{\mu_{1,j}}$. Furthermore, in conjunction with \autoref{cor:description-factor}, it is straightforward to compute the dimension of $M_i / M_{i+1}$ for any $1 \leq i \leq r - 1$.
\bigskip
The remainder of this section is devoted to proving the auxiliary results used in the proof of \autoref{thm:specht-filtration-general-partition}. For convenience, we adopt the notation and assumptions of \autoref{thm:specht-filtration-general-partition} throughout the remainder of this section, unless stated otherwise. We recall and fix some of them here for clarity. For $\lambda=(\lambda_1,\cdots,\lambda_r)\in\Par[\Lambda]_{\alpha}$, we set $\beta=\alpha_{(\lambda_1)}\in Q^+$ and $\nu=(\lambda_2, \dots, \lambda_r)$. The module $S^\beta:=S^{(\lambda_1)}$ is the one-dimensional Specht module associated with $(\lambda_1)$ over $R^{\Lambda}_\beta$. Assume $M^{\nu}$ has a generalized Specht filtration as in \autoref{eq:generalized-filtration} where $N_i$ is the $i$-th stage. 
%%%%%%%%%%

%%%%%%%%%%
\begin{Lemma}\label{lm:from-Ni-to-Mi+1}
    For $1\leq i\leq r-2$, we have $M_{i+1}\cong S^\beta\circ N_i$ and $M^\lambda=S^\beta\circ M^\nu$.
\end{Lemma}
\begin{proof}
The last statement follows directly from the definition. Recall that $N_i$ is an $R_{\alpha - \beta}$-submodule of $M^\nu$.

For $1 \leq i \leq r - 2$, let $B_{i+1} = (i+1, \lambda_{i+2}) \in [\lambda]$ and $C_i = (i, \lambda_{i+2}) \in [\nu]$ denote the last Garnir nodes in the $(i+1)$-st row of $[\lambda]$ and the $i$-th row of $[\nu]$, respectively. Then, by definition, we have
\[
    M_{i+1} = R_\alpha \{ \psi^{B_{i+1}} v, \dots, \psi^{B_{r-1}} v \} \subset M^\lambda
\]
and
\[
    N_i = R_{\alpha - \beta} \{ \psi^{C_i} v, \dots, \psi^{C_{r-2}} v \} \subset M^\nu.
\]
Let $v_\beta$ and $v_N$ denote the standard cyclic generators of $L_\beta = S^\beta$ and $M^\nu$, respectively. Let $v_M := v_\beta \otimes v_N$ be the standard cyclic generator of $M^\lambda$.

Applying the exact functor $L_\beta \circ - := \operatorname{Ind}^{R_\alpha}_{R_{\beta, \alpha - \beta}} (L_\beta \otimes -)$ to $N_i \subset M^\nu$, we obtain
\[
    L_\beta \circ N_i \subset L_\beta \circ M^\nu \cong M^\lambda,
\]
where the isomorphism on the right maps $v_\beta \otimes v_N$ to $v_M$.

By the standard inclusion $R_\beta \otimes R_{\alpha - \beta} \hookrightarrow R_\alpha$, we know that
\[
    v_\beta \otimes \psi^{C_i} v_N = \psi^{B_{i+1}} (v_\beta \otimes v_N)
\]
for each $1 \leq i \leq r - 2$. Hence,
\begin{align*}
    L_\beta \circ N_i 
    &= R_\alpha \{ v_\beta \otimes \psi^{C_i} v_N, \dots, v_\beta \otimes \psi^{C_{r-2}} v_N \} \\
    &= R_\alpha \{ \psi^{B_{i+1}} (v_\beta \otimes v_N), \dots, \psi^{B_{r-1}} (v_\beta \otimes v_N) \} \\
    &\cong R_\alpha \{ \psi^{B_{i+1}} v_M, \dots, \psi^{B_{r-1}} v_M \} \\
    &= M_{i+1},
\end{align*}
as desired.
\end{proof} 
%%%%%%%%%%
For simplicity, from now on, we simplify the notations and write $k:=k_1$, $\mu_j:=\mu_{1,j}$ and $\upphi_j:=\upphi_{1,j}$ where $1\leq j\leq k$.

Recall that for a cyclotomic KLR algebra $R^\Lambda_\alpha$, where $\Lambda$ is of level $\ell$, we can associate to each $\ell$-partition $\blam \in \Par_\alpha$ a cell module $C^\blam$, which is (graded) isomorphic to the Specht module $S^\blam$. There exists a distinguished subset of $\Par_\alpha$ called the set of \emph{Kleshchev partitions}. If $\blam$ is a Kleshchev partition, then $S^\blam$ is indecomposable and $D^\blam$ is its unique irreducible head. Moreover, the set 
\[
\{D^\blam \mid \blam \in \Par_\alpha \text{ is Kleshchev}\}
\]
is a complete set of irreducible modules for $R^\Lambda_\alpha$. These results can be found in \cite{humathas-graded-cellular}.

For our purposes, we do not need the recursive definition of Kleshchev partitions; instead, we record the following results.
\begin{Proposition}[{\cite[Corollary 3.23]{humathas-quiver-schur-linear-quiver}}]\label{prop:kleshchev-partition}
    Suppose that $e = 0$ or $e > n$, $\kappa_1 \geq \kappa_2 \geq \dots \geq \kappa_\ell$ and $\mu \in \Par[\Lambda]_n$. 
    Then $\mu = (\mu^{(1)}, \dots, \mu^{(\ell)})$ is Kleshchev if and only if
    \[
        \mu^{(l)}_{r + \kappa_l - \kappa_{l+1}} \leq \mu^{(l+1)}_{r}, \quad \text{for } 1 \leq l < \ell \text{ and } r \geq 1.
    \]
\end{Proposition}
%%%%%%%%
\begin{Definition}
    Let $\Std^\bmu(\blam)$ be the set
    $\{\bs\in\Std(\blam)|\bs\unrhd T^\bmu \text{ and }\res(\bs)=\bi^\bmu\}$.    
\end{Definition}
\begin{Definition}\label{dfn:graded-decomposition-number}
    Suppose that $\blam, \bmu\in\Par_{n}$. Define the graded decomposition number to be
  $$d_{\blam\bmu}(q)=[S^\blam:D^\bmu]_q
                    =\sum_{d\in\Z}\, [S^\blam:D^\bmu\<d\>]\,q^d.$$
    where $[M:L]$ is the graded multiplicity of $L$ in $M$ for any graded simple module $L$ and graded module $M$.
\end{Definition}
\begin{Proposition}[{\cite[Appendix B]{humathas-quiver-schur-linear-quiver}}]\label{prop:decomposition-number-level-2}
     Fix a linear quiver $A_\infty$ or $\Aone[e-1]$ with $e\gg 0$. 
    Suppose $\Lambda$ is of level $2$, then $\#\Std^\bmu(\blam)\leq 1$. If the equality holds, let $t^\bmu_\blam$ be the unique element in $\Std^\bmu(\blam)$.
    Moreover, (suppose $\bmu$ is a Kleshchev partition), we have:
    \[
        d_{\blam\bmu}(q)=\begin{cases}q^{\deg t^\bmu_\blam-\deg T^\bmu} & \text{ if }\#\Std^\bmu(\blam)=1\\ 0& \text{else}\end{cases}
    \]  
\end{Proposition}
%%%%%%%%%
\begin{Lemma}\label{lm:mu_j-is-kleshchev-partitions}
    For $1\leq j\leq k_1$, $\mu_j$ is a Kleshchev partition.
\end{Lemma}
\begin{proof}
    By \autoref{prop:kleshchev-partition}, we only need to verify ${(\mu_j)}_{r+1}^{(1)}\leq {(\mu_j)}_{r}^{(2)}$ for each $r\geq 1$. By construction in \autoref{thm:specht-filtration-general-partition}, we know this is true since ${\mu_j}^{(1)}$ consists of one row.
\end{proof}
%%%%%%%%%%%
We introduce two useful quantities for a partition $\lambda=(\lambda_1,\cdots,\lambda_r)$. For all admissible $i$, set
\begin{equation}\label{eq:ni-di}
    n_i:=\sum\limits_{1\leq j\leq i}\lambda_j,\qquad d_i=\lambda_i-\lambda_{i+1}+1.
\end{equation}
%%%%%%%%%%
\begin{Lemma}\label{lm:dominant-resolution}
    The Specht resolution in \autoref{eq:specht-resolution-thm} when $i=1$ are just the $2$-partitions $\mu$ in $\Par[\Lambda(1)]_\alpha$ such that $\mu\trianglelefteq \mu_1$ listed in the dominance order of partitions, i.e. $\mu_{j}\triangleleft
\mu_{i}$ if and only if $j> i$. 
\end{Lemma}
\begin{proof}
    Since we are working in type $A_\infty$ or $\Aone[e-1]$ with $e \gg 0$, each diagonal of $[\lambda]$ has a distinct residue. The only possible way to move a removable node up or down while preserving the residue is to move it along a diagonal, thereby keeping the content unchanged. It is then easy to see that, in order to maintain the structure of a $2$-partition, it is impossible to move any node within a single component. 

    Hence, the only way to construct $2$-partitions strictly smaller than $\mu_1 = (\lambda_2 - 1 \mid \nu_1)$—where $\nu_1 = (\lambda_1 + 1, \lambda_3, \cdots, \lambda_r)$—is to move nodes from the first component to the second. Under our chosen residue sequence, the only such possibility is to move the last $d_1$ nodes into the second row of $\nu_1$, which yields $\mu_2$.

    Similarly, for each $j$, the partition $\mu_{j+1}$ is the unique $2$-partition in $\Par[\Lambda(i)]_\alpha$ that lies immediately below $\mu_j$ with respect to the dominance order. The number $k$ is the maximal index such that there does not exist any $2$-partition in $\Par[\Lambda(i)]_\alpha$ lying strictly below $\mu_k$.
\end{proof}
%%%%%%%%%%
\begin{Corollary}\label{cor:decomposition-number-dominance}
    For $\mu_j,1\leq j\leq k$, we have:
    \[
    d_{\mu_{j}\nu}(q)=
    \begin{cases}
        q^{\deg t^\nu_{\mu_j}-\deg T^\nu} & \text{ if }\nu=\mu_{j}\text{ or }\mu_{j+1}\\ 
        0& \text{else}
    \end{cases}
    \]
\end{Corollary}
\begin{proof}
    By \autoref{prop:decomposition-number-level-2}, we only need to show that $\#\Std^\nu(\mu_j) = 1$ if and only if $\nu = \mu_j$ or $\nu = \mu_{j+1}$. 

    If $\nu = \mu_j$, this is clear by taking $t^\nu_{\mu_j} = T^{\mu_j}$. 

    If $\nu = \mu_{j+1}$, set $n_j$ and $d_j$ as in \autoref{eq:ni-di}. It is straightforward to verify that $t^\nu_{\mu_j}$ is obtained from $T^{\mu_{j+1}}$ by moving the last $d_j$ nodes from the $(j+1)$-st row to the first component, concatenating them with the first row. In other words, $t^\nu_{\mu_j}$ is the following tableau:
    
        \Tableau[xscale=4.2,dotted cols={2,5}]{1\dots{\lambda_{j+2}{-}j{-}1}*{2\lambda_{j+2}{+}n_j}*\dots*{n_{j+2}}}\\
        %%%%%%%%%%%%%%%%%%%%%%%%%%%%%%%%%
        \hspace*{3mm}\Tableau[xscale=4.2,dotted cols={3,5},dotted rows={2,3,6,7}]{{\lambda_{j+2}{-}j}\dots 0\dots 0\dots\dots{\lambda_{j+2}{-}j{+}\lambda_1},0000000,0000000,{\lambda_{j+2}{+}n_{j-1}{-}1}\dots 0\dots 0 {\lambda_{j+2}{+}n_j{-}1},{\lambda_{j+2}{+}n_j}\dots 0{2\lambda_{j+2}{+}n_j{-}1},000,00,\dots}
    
    By the standard theory of cellular algebras (see \cite{humathas-graded-cellular}, for example), we have $d_{\mu_j,\nu} \neq 0$ only if $\mu_j \trianglerighteq \nu$. Hence, by \autoref{lm:dominant-resolution}, it suffices to verify that there is no element in $\Std^\nu(\mu_j)$ for $\nu = \mu_s$ with $s \geq j+2$.

    Suppose, for contradiction, that there exists $T \in \Std^\nu(\mu_j)$. Then the first $\lambda_{s+1} - s$ entries must be $1, 2, \dots, \lambda_{s+1} - s$. Furthermore, by the condition $\bi^T = \bi^{T^\nu}$, the first $j$ rows of the second component must coincide with $T^\nu$ as well: entries increase (in row-reading order) from $\lambda_{s+1}-s+1$ to $\lambda_{s+1-s+n_j+j}$. The only possible difference begins at the last $d_{j+1}$ nodes of the $(j+1)$-st row of $T^\nu$: in the second component, $\mu_{j}$ has $\lambda_{j+2}$ nodes in the $(j+1)$-st row, whereas $\nu$ has $\lambda_{j+2} + d_{j+1}$ nodes. In this row of $T$, the first $\lambda_{j+2}$ entries are identical to those in $T^\nu$. However, to satisfy $T \trianglerighteq T^\nu$, the next node must lie in the first component. But the node with this residue is not adjacent to the $(\lambda_{s+1} - s)$-th node unless $s = j+1$. Therefore, the desired tableau cannot be standard.
\end{proof}
\begin{Corollary}\label{cor:minimal-element}
    $S^{\mu_k}$ is irreducible.
\end{Corollary}
\begin{proof}
    Since $\mu_k$ is a minimal element in $\Par[\Lambda(1)]_\alpha$, this can be verified either using \autoref{cor:decomposition-number-dominance}, or deduced from the standard theory of cellular algebras together with \autoref{lm:dominant-resolution}.
\end{proof}
%%%%%%%%%%
%%%%%%%%%
\begin{Lemma}\label{lm:kernel-not-zero}
If $k\neq 1$, then $\psi^{t^{\mu_{2}}_{\mu_{1}}}\psi^{B_1}v\in M_2$.
\end{Lemma}
\begin{proof}
    The condition $k\neq 1$ is equivalent to $\lambda_3\geq 2$. Set $n_i$ and $d_i$ as in \autoref{eq:ni-di}. 
    Let $T:=t^{\mu_{2}}_{\mu_{1}}$, by \autoref{cor:decomposition-number-dominance}, we know $T$ is of the following form:
    %%%%%
    
        \Tableau[xscale=3.3,dotted cols={2,5}]{10{\lambda_3{-}2}{\lambda_1{+}2\lambda_3}0{n_3}}\\
        \hspace*{3mm}\Tableau[xscale=3.3,dotted cols={3,4,5,7,8},dotted rows={3,4,5}]{{\lambda_3{-}1}{\lambda_3}000{2\lambda_3{-}2}00{\lambda_3{+}\lambda_1{-}1},{\lambda_3{+}\lambda_1}{\lambda_3{+}\lambda_1{+}1}000{\lambda_1{+}2\lambda_3{-}1},000,00,00}\\

    %\Multitableau[xscale=2.8,dotted cols={2,5|3,5},dotted rows={|3,4,5}]{10{\lambda_3{-}2}{\lambda_1{+}2\lambda_3}0{n_3}|{\lambda_3{-}1}{\lambda_3}0{2\lambda_3{-}2}0{\lambda_3{+}\lambda_1{-}1},{\lambda_3{+}\lambda_1}{\lambda_3{+}\lambda_1{+}1}0{\lambda_1{+}2\lambda_3{-}1},0,0,0}
    
    Hence (for simplicity, we omit the idempotent in the expression of $\psi^T$)
    \[
        \psi^T=(\psi_{n_3-d_2}\cdots\psi_{\lambda_3-1})\cdots (\psi_{n_3-1}\cdots\psi_{\lambda_2-1}).
    \]
    The Garnir tableau of $B_1$ is of the following form:
    \begin{center}
        \resizebox{0.9\textwidth}{!}{\Tableau[xscale=2.5,dotted cols={3,4,8,9},dotted rows={3,4,5}]{1200{\lambda_2-1}*{2\lambda_2}*{2\lambda_2{+}1}00*{\lambda_1+\lambda_2},*{\lambda_2}*{\lambda_2{+}1}00*{2\lambda_2{-}2}*{2\lambda_2{-}1},000,00,00}}
    \end{center}
    and
    \[
        \psi^{B_1}=(\psi_{2\lambda_2-1}\cdots\psi_{n_2-1})\cdots(\psi_{\lambda_{2}}\cdots\psi_{\lambda_1})e(\bi).
    \]
    Set $A_i=(1,i)\in [\lambda]$ for $1\leq i\leq \lambda_2$, then $A_i$ are all the Garnir nodes between the first two rows and $B_1=A_{\lambda_2}$. We will keep applying \autoref{lm:garnir-connection-2} (modulo $\pm$ signs):
    \begin{align*}
        &\psi^T\psi^{B_1}v\\
        %%%%
        &=(\psi_{n_3-d_2}\cdots\psi_{\lambda_3-1})\cdots(\psi_{n_3-1}\cdots\psi_{\lambda_2-1})\psi^{B_1}v\\[2pt]
        %%%%
        &=(\psi_{n_3-d_2}\cdots\psi_{\lambda_3-1})\cdots(\psi_{n_3-1}\cdots\psi_{\lambda_2-1})\psi^{A_{\lambda_2}}v\\[2pt]
        %%%%
        &=(\psi_{n_3-d_2}\cdots\psi_{\lambda_3-1})\cdots
          (\psi_{n_3-1}\cdots\psi_{\lambda_1+\lambda_2})
          (\psi_{\lambda_1+\lambda_2-1}\cdots\psi_{\lambda_2-1})
          \psi^{A_{\lambda_2}}v\\[2pt]
            %%%%
        &=(\psi_{n_3-d_2}\cdots\psi_{\lambda_3-1})
          \cdots
          (\psi_{n_3-1}\cdots\psi_{\lambda_1+\lambda_2})
          \psi^{A_{\lambda_2-1}}v\\[2pt]
            %%%%
        &=(\psi_{n_3-d_2}\cdots\psi_{\lambda_3-1})
          \cdots(\psi_{n_3-2}\cdots\psi_{\lambda_2-2})
          (\psi_{n_3-1}\cdots\psi_{\lambda_1+\lambda_2})
          \psi^{A_{\lambda_2-1}}v\\[2pt]
            %%%%
        &=(\psi_{n_3-d_2}\cdots\psi_{\lambda_3-1})
          \cdots(\psi_{n_3-2}\cdots\psi_{\lambda_1+\lambda_2-1})(\psi_{\lambda_1+\lambda_2-2}\cdots\psi_{\lambda_2-2})
          (\psi_{n_3-1}\cdots\psi_{\lambda_1+\lambda_2})
          \psi^{A_{\lambda_2-1}}v\\[2pt]
            %%%%
        &=(\psi_{n_3-d_2}\cdots\psi_{\lambda_3-1})
          \cdots(\psi_{n_3-2}\cdots\psi_{\lambda_1+\lambda_2-1})(\psi_{n_3-1}\cdots\psi_{\lambda_1+\lambda_2})(\psi_{\lambda_1+\lambda_2-2}\cdots\psi_{\lambda_2-2})
          \psi^{A_{\lambda_2-1}}v\\[2pt]
            %%%%
        &=(\psi_{n_3-d_2}\cdots\psi_{\lambda_3-1})
          \cdots(\psi_{n_3-2}\cdots\psi_{\lambda_1+\lambda_2-1})(\psi_{n_3-1}\cdots\psi_{\lambda_1+\lambda_2})
          \psi^{A_{\lambda_2-2}}v\\[2pt]
            %%%%%
        &=\cdots\\
            %%%%%
        &=(\psi_{n_3-d_2}\cdots\psi_{\lambda_3+\lambda_1})\cdots(\psi_{n_3-1}\cdots\psi_{\lambda_1+\lambda_2})\psi^{A_{\lambda_3-1}}v\\
        %%%%
        &=(\psi_{n_3-d_2}\cdots\psi_{n_3-1})\cdots(\psi_{\lambda_3+\lambda_1}\cdots\psi_{\lambda_1+\lambda_2})\psi^{A_{\lambda_3-1}}v\\
        %%%%
        %%%%
        &=\psi^{B_2}\psi^{A_{\lambda_3-1}}v\\
        &=\psi^{A_{\lambda_3-1}}\psi^{B_2}v\in M_2.
    \end{align*}
     Note $B_2=(2,\lambda_3)\in [\lambda]$ and $\psi^{B_2}$ commutes with $\psi^{A_{\lambda_3}-1}$ because the two Garnir belts do not intersect and hence all the $\psi_i$ in the two expressions commute.
\end{proof}
%%%%%%%%%%%%%%%%%%%%%%%%%%%%%%%%%%%%%%%%%%%%%%%%%%%%%%%%%%%%%%%%%%%%%%%%%%%%%%%%%%%%%%%%%%%%%%%%%%%%%%%%%%%%%%%%%%%%%%%%%%%%%%%%%%%%%%%%%%%%%%%%%%%%%%%%%%%%%%%%%%%%%%%%%%%%%%%%%%%%%%%%%%%%%%%%%%%%%%%%%%%%%%%%%%%%%%%%%%%%%%%%%%%%%%%%%%%%%%%%%%%%%%%%%%%%%%%%%%%%%%%%%%%%%%%%%%%%%%%%%%%%%%%%%%%%%%%%%%%%%%%%%%%%%%%%%%%%%%%%%%%%%%%%%%%%%%%%%%%%%%%%%%%%%%%%%%%%%%%%%%%%%%%%%%%%%%%%%%%%%%%%%%%%%%%%%%%%%%%%%%%%%%%%%%%%%%%%%%%%%%%%%%%%%%%%%%%%%%%%%%%%%%%%%%%%%%%%%%%%%%%%%%%%%%%%%%%%%%%%%%%%%%%%%%%%%%%%%%%%%%%%%%%%%%%%%%%%%%%%%%%%%%%%%%%%%%%%%%%%%%%%%%%%%%%%%%%%%%%%%%%%%%%%%%%%%%%%%%%%%%%%%%%%%%%%%%%%%%%%%%%%%%%%%%%%%%%%%%%%%%%%%%%%%%%%%%%%%%%%%%%%%%%%%%%%%%%%%%%%%%
\section{Higher Levels and Skew Specht Filtrations}\label{sec:higher-level-skew}
%%%%%%%%%%%%%%%%%%%%%%
\subsection{Higher Level Case}

In this section, we briefly discuss how to construct a (generalized) Specht filtration of $M^\blam$ for $\blam$ an $\ell$-partition with $\ell > 1$. 

Let $\blam \in \Par_\alpha$, where $\Lambda$ is a dominant weight of level $\ell$. By definition, we have 
\[
    M^\blam \cong M^{\blam^{(1)}} \circ \cdots \circ M^{\blam^{(\ell)}}.
\]
For each $1 \leq s \leq \ell$, we have constructed a (generalized) Specht filtration of $M^{\blam^{(s)}}$ in various cases, as described in \autoref{sec:filtration}, \autoref{sec:specht-filtration-2row}, and \autoref{sec:specht-filtration-general-partition}. Suppose the length of $\blam^{(s)}$ is $r_s$, with the filtration of $M^{\blam^{(s)}}$ given by:
\[
    M^{\blam^{(s)}} = M^{s}_{0} \supsetneq M^{s}_{1} \supsetneq \cdots \supsetneq M^{s}_{r_s-1} \supsetneq M^{s}_{r_s} = 0.
\]
For each $1 \leq i \leq r_s$, there exists a Specht resolution:
\[
    0 \to S^{\mu^s_{i k_{si}}} \to \cdots \to S^{\mu^s_{i1}} \twoheadrightarrow M^s_i / M^s_{i+1} \to 0,
\]
where $k_{si} = k_i(\blam^{(s)})$.

Define
\[
    M_i := M^1_i \circ M^2_i \circ \cdots \circ M^\ell_i,
\]
so that we obtain a filtration of $M^\blam$:
\[
    M^\blam = M_0 \supsetneq M_1 \supsetneq \cdots \supsetneq M_{r-1} \supsetneq M_r = 0,
\]
where $r = \max\{r_s \mid 1 \leq s \leq \ell\}$, and we set $M^s_i := \bk$ to be the trivial module whenever $i > r_s$.

Since the external tensor product commutes with taking quotients, we have:
\[
    M_i / M_{i+1} \cong M^1_i / M^1_{i+1} \circ M^2_i / M^2_{i+1} \circ \cdots \circ M^\ell_i / M^\ell_{i+1}.
\]
Hence, for each $i$, we obtain the following Specht resolution:
\[
    0 \to S^{\mu^1_{i k_i}} \circ S^{\mu^2_{i k_i}} \circ \cdots \circ S^{\mu^\ell_{i k_i}} \to \cdots \to 
    S^{\mu^1_{i1}} \circ S^{\mu^2_{i1}} \circ \cdots \circ S^{\mu^\ell_{i1}} \twoheadrightarrow M_i / M_{i+1} \to 0,
\]
where $k_i = \max\{k_{si} \mid 1 \leq s \leq \ell\}$. For any $j > k_{si}$, we set $\mu^s_{ij} := \emptyset$ to be the empty partition, so that $S^{\mu^s_{ij}}$ is the trivial module.

For each $1 \leq j \leq k_i$, define
\[
    \nu_{ij} := (\mu^1_{ij} \mid \mu^2_{ij} \mid \cdots \mid \mu^\ell_{ij}) \in \Par_\alpha.
\]
Then, by \autoref{thm:kmr-high-level-specht}, the above Specht resolution becomes
\[
    0 \to S^{\nu_{i k_i}} \to \cdots \to S^{\nu_{i1}} \twoheadrightarrow M_i / M_{i+1} \to 0.
\]
%%%%%%%%%%%%%%%%%%%%%%
\subsection{Skew Specht filtration}
As mentioned in \autoref{rmk:skew-specht-filtration-hook}, we can construct a skew Specht filtration by reversing the order of the Garnir relations. To be precise, suppose we are working in type~$A_\infty$, and let $\lambda = (\lambda_1, \dots, \lambda_r) \in \Par_\alpha$ be a partition. Define $C_i := (i,1) \in [\lambda]$ for $1 \leq i \leq r-1$ to be the set of Garnir nodes in the first column. By \autoref{thm:2row-garnir-connection}, we know that 
\[
M_1 = R_\alpha \{ \psi^{C_1}v, \dots, \psi^{C_{r-1}}v \}
\]
is the submodule generated by all Garnir relations. Let $v$ be the standard cyclic generator of $M^\lambda$. Instead of defining 
\[
M_i := R_\alpha \{\psi^{C_i}v, \dots, \psi^{C_{r-1}}v \} \quad \text{for } 1 \leq i \leq r-1,
\]
as in \autoref{thm:specht-filtration-general-partition}, we define
\[
M'_i := R_\alpha \{\psi^{C_1}v, \dots, \psi^{C_{r-i}}v \}.
\]
Then the filtration
\[
M^\lambda = M'_0 \supsetneq M'_1 \supsetneq \cdots \supsetneq M'_{r-1} \supsetneq M'_r = 0
\]
has the following property: for each $1 \leq i \leq r-1$, there exists a resolution of $M'_i/M'_{i+1}$ by skew Specht modules as introduced in \cite{muth-graded-skew-specht}. The explicit formulas for the skew partitions that appear in this resolution, as well as the detailed proof, are omitted. 
\bibliography{reference}  % Without the .bib extension
\bibliographystyle{alpha}  % Or another style of your choice
\end{document}